\documentclass[12pt,letterpaper]{amsart}
\usepackage[T1]{fontenc}
\usepackage[utf8]{inputenc}
\usepackage{lmodern}
\usepackage[letterpaper,margin=1in,headheight=14pt,headsep=15pt,footskip=14pt]{geometry}
\usepackage{amssymb,mathtools}
\usepackage{booktabs}
\usepackage{microtype}
\usepackage{hyperref}
\hypersetup{
            colorlinks=true,
            linkcolor=red,
            anchorcolor=blue,
            citecolor=green}

\numberwithin{equation}{section}
\newtheorem{theorem}{Theorem}[section]
\newtheorem{lemma}[theorem]{Lemma}
\newtheorem{proposition}[theorem]{Proposition}
\newtheorem{corollary}[theorem]{Corollary}
\theoremstyle{remark}

\newcommand{\Sym}{\operatorname{Sym}}
\newcommand{\SL}{\operatorname{SL}}
\newcommand{\CM}{\mathrm{CM}}
\newcommand{\ST}{\mathrm{ST}}
\newcommand{\Res}{\operatorname*{Res}}
\newcommand{\C}{\mathbf C}
\newcommand{\R}{\mathbf R}

\newcommand{\e}{\mathrm e}

\allowdisplaybreaks[1]
\title[Moments over Friable Integers]
{Moments of Symmetric-Power Coefficients
at Powers of Friable Integers}
\author{Nanxiang Wang}
\subjclass[2020]{Primary 11N37;
Secondary 11N25, 11F30, 11F66.}
\keywords{Friable integers, Hecke eigenvalues, symmetric powers, moments, Dickman function, complex multiplication}
\thanks{N. X. Wang, School of Mathematics and Statistics, Qingdao University, 308 Ningxia Road, Qingdao, Shandong 266071, P. C. China, \href{mailto:wangnanxiang@qdu.edu.cn}{\texttt{wangnanxiang@qdu.edu.cn}}.}
\date{}
\begin{document}

\begin{abstract}
Let \(f\) be a fixed normalized holomorphic cuspidal Hecke eigenform of level one.
For fixed positive integers \(j,k,\ell\), we study the sum of
\(\lambda_{\Sym^j f}(n^k)^\ell\) over the \(y\)-friable integers \(n\le x\),
that is, integers whose prime factors do not exceed \(y\).
The prime-mean parameter is the multiplicity of the trivial representation in
\((\Sym^k(\Sym^j\C^2))^{\otimes\ell}\), which gives an exact procedure for computing
the logarithmic exponent. For even \(\ell\), we obtain an asymptotic formula for
\(y\ge\exp\{A(\log\log x)^2\}\). For signed moments with positive prime mean, we
obtain an additive main-term formula for
\(y\ge\exp\{A(\log x\log\log x)^{2/3}\}\), where \(A\) is sufficiently large and
depends on the fixed data. Two explicit odd moments have positive constants and
therefore have a matching upper and lower order. For zero prime mean, an
exponential-convolution argument gives exponential cancellation. We also treat
ordinary moments of fixed holomorphic newforms with complex multiplication
(CM): the even-moment parameter is
\(\frac12\binom{2r}{r}\), and a Hecke zero-free region gives different friable
ranges. 
\end{abstract}

\maketitle

\section{Introduction and main results}\label{sec:intro}

An integer \(n\ge1\) is called \(y\)-friable if its largest prime factor
\(P^+(n)\) does not exceed \(y\); throughout, \(P^+(1)=1\). We write
\[
 \Psi(x,y)=\#\{n\le x:P^+(n)\le y\},\qquad
 X=\log x,\quad L=\log y,\quad u=\frac{X}{L}.
\]
All logarithms are natural. The asymptotic statements below concern
\(x\to\infty\) with \(2\le y\le x\). Constants denoted by \(A\) in different
statements need not coincide; each is sufficiently large and depends only on
the fixed parameters specified there.

Friable mean values of multiplicative functions depend on their mean at primes
and on their local factors at higher prime powers. The general nonnegative
theory of Tenenbaum and Wu~\cite{TW}, the analytic evaluation of de Bruijn main
terms by Hanrot, Tenenbaum and Wu~\cite{HTW}, and the complex-valued theory of
de la Bretèche and Tenenbaum~\cite{BT} provide the relevant background. We
apply these ideas to the combined family
\begin{equation}\label{eq:sum}
 S_{j,k,\ell}(x,y;f)
   =\sum_{\substack{n\le x\\P^+(n)\le y}}
      \lambda_{\Sym^j f}(n^k)^\ell .
\end{equation}

\subsection{Notation}

Unless stated otherwise, \(f\) is a fixed normalized holomorphic cuspidal Hecke
eigenform of even weight \(w\ge2\) for \(\SL_2(\mathbf Z)\). In particular, \(f\)
is non-CM. We use the normalization
\[
 f(z)=\sum_{n\ge1}\lambda_f(n)n^{(w-1)/2}\e^{2\pi inz},
 \qquad \lambda_f(1)=1.
\]
The symbol \(\lambda_{\Sym^j f}(n)\) denotes the normalized Dirichlet coefficient
of the standard \(L\)-function \(L(s,\Sym^j f)\). The positive integers
\(j,k,\ell\) are fixed; no uniformity in these integers or in \(f\) is asserted.
All these coefficients are real under the present hypotheses.
The letter \(p\) denotes a prime. We write \(\mathbf 1_E\) for the
indicator of a set \(E\). In statements for a general arithmetic function,
implied constants may depend on the fixed data in that statement; no
additional uniformity is implicit.

At a prime \(p\), choose \(\theta_p\in[0,\pi]\) so that the normalized Satake
parameters are \(\alpha_p=\e^{i\theta_p}\) and \(\alpha_p^{-1}\).
Deligne's theorem~\cite{Deligne} gives \(|\alpha_p|=1\), and
\begin{equation}\label{eq:local}
 \sum_{\nu\ge0}\lambda_{\Sym^j f}(p^\nu)T^\nu
   =\prod_{r=0}^{j}(1-\alpha_p^{\,j-2r}T)^{-1}.
\end{equation}
We set \(L(s,\Sym^0f)=\zeta(s)\). Let \(U_m\) be the Chebyshev polynomial of the
second kind, normalized by
\[
 U_0(t)=1,\quad U_1(t)=2t,\quad
 U_{m+1}(t)=2tU_m(t)-U_{m-1}(t).
\]
Equivalently,
\[
 U_m(\cos\theta)=\frac{\sin((m+1)\theta)}{\sin\theta}
 \quad(0<\theta<\pi),
\]
with continuous endpoint values. Thus
\(\lambda_{\Sym^m f}(p)=U_m(\cos\theta_p)\).
For \(z>0\), the generalized divisor function \(d_z\) is defined by
\[
 \zeta(s)^z=\sum_{n\ge1}\frac{d_z(n)}{n^s}\quad(\Re s>1),
 \qquad
 d_z(p^\nu)=\frac{\Gamma(z+\nu)}{\Gamma(z)\nu!}.
\]
Here \(\zeta(s)^z=\exp\{z\log\zeta(s)\}\), using the logarithm given
by the absolutely convergent Euler product on \(\Re s>1\).
For an integer \(D\ge1\), this gives
\(d_D(p^\nu)=\binom{D+\nu-1}{\nu}\).

Define the Laurent polynomial
\begin{equation}\label{eq:Q}
 Q_{j,k}(z)=
 \sum_{\substack{a_0+\cdots+a_j=k\\a_r\in\mathbf Z_{\ge0}}}
    z^{\,jk-2\sum_{r=0}^{j}ra_r}.
\end{equation}
Here and below, \([z^r]\) denotes coefficient extraction. Put
\begin{equation}\label{eq:index}
 \kappa=\kappa_{j,k,\ell}
   =[z^0]Q_{j,k}(z)^\ell-[z^2]Q_{j,k}(z)^\ell,
 \qquad
 D=\binom{j+k}{j}^{\ell}.
\end{equation}
Section~\ref{sec:representation} identifies \(\kappa\) as an invariant
multiplicity; in particular, it is an integer with \(0\le\kappa\le D\).

For \(\beta>0\), let \(\rho_\beta\) be the continuous, locally absolutely
continuous function on \((0,\infty)\) specified successively on unit intervals by
\begin{equation}\label{eq:rho}
 \begin{split}
 \rho_\beta(v)&=\frac{v^{\beta-1}}{\Gamma(\beta)}
       \qquad(0<v\le1),\\
 v\rho_\beta'(v)+(1-\beta)\rho_\beta(v)
       +\beta\rho_\beta(v-1)&=0
       \qquad\text{for almost every }v>1,
 \end{split}
\end{equation}
with zero extension to \(v<0\). We take \(\rho_1(0)=1\) and
\(\rho_\beta(0)=0\) for \(\beta>1\). When \(0<\beta<1\), the function
has an integrable singularity at zero; its assigned value at zero is
irrelevant to the Lebesgue integrals below and may be set to zero.
We write \(\rho=\rho_1\). The almost-everywhere convention avoids imposing
a derivative at a junction where the one-sided derivatives need not agree.
Basic positivity and convolution properties are proved in
Lemma~\ref{lem:dickman-basic}.

We distinguish two convolution operations by their arguments. For arithmetic
functions,
\[
 (a*b)(n)=\sum_{d\mid n}a(d)b(n/d)
\]
is Dirichlet convolution. For functions on \([0,\infty)\), extended by zero
to negative arguments, additive convolution is
\[
 (g*h)(v)=\int_0^v g(v-t)h(t)\,dt.
\]
For locally finite signed measures, additive convolution is the pushforward
of their product measure under addition, whenever the integral is defined;
in particular,
\[
 (\mu*g)(v)=\int_{[0,v]}g(v-t)\,d\mu(t).
\]
The notation \(g^{*r}\) denotes the \(r\)-fold additive convolution when
\(g\) is a function of a real variable. Our Laplace-transform convention is
\[
 \widehat g(s)=\int_0^\infty \e^{-sv}g(v)\,dv,
 \qquad
 \widehat\mu(s)=\int_{[0,\infty)}\e^{-sv}\,d\mu(v),
\]
on any half-plane of absolute convergence. The second integral includes
any atom at zero. A right-continuous function \(g\) of locally bounded
variation, zero on \(( -\infty,0)\), defines the Stieltjes measure \(dg\)
by \(dg((a,b])=g(b)-g(a)\).

\subsection{Main results}

\begin{theorem}[Even moments]\label{thm:even}
Suppose that \(\ell\) is even. Then \(\kappa\ge1\), and there is a constant
\(A=A(f,j,k,\ell)>0\) such that, uniformly as \(x\to\infty\) in
\begin{equation}\label{eq:even-range}
 \exp\{A(\log\log x)^2\}\le y\le x,
\end{equation}
one has
\begin{equation}\label{eq:even}
 S_{j,k,\ell}(x,y;f)
  =C_{j,k,\ell}(f)xL^{\kappa-1}\rho_\kappa(u)
    \left\{1+O_{f,j,k,\ell}
       \left(\frac{\log(u+1)}{L}\right)\right\}.
\end{equation}
The constant is strictly positive and is given by
\begin{equation}\label{eq:constant}
 C_{j,k,\ell}(f)
  =\lim_{z\to\infty}\prod_{p\le z}
       (1-p^{-1})^\kappa
       \left(\sum_{\nu\ge0}
          \frac{\lambda_{\Sym^j f}(p^{k\nu})^\ell}{p^\nu}\right).
\end{equation}
The limit is taken in increasing prime order.
\end{theorem}

The product in \eqref{eq:constant} is generally not absolutely convergent.
Lemma~\ref{lem:euler} proves its boundary convergence and identifies it with an
analytic value; an absolutely convergent correction product appears there.

\begin{theorem}[Positive prime mean]\label{thm:signed}
Suppose that \(\kappa\ge1\). There are a real constant
\(C=C_{j,k,\ell}(f)\), given by \eqref{eq:constant}, and a constant
\(A=A(f,j,k,\ell)>0\) such that, uniformly in
\begin{equation}\label{eq:signed-range}
 \exp\{A(\log x\log\log x)^{2/3}\}\le y\le x,
\end{equation}
one has the additive formula
\begin{equation}\label{eq:signed}
 S_{j,k,\ell}(x,y;f)
  =xL^{\kappa-1}\rho_\kappa(u)
      \left\{C+O_{f,j,k,\ell}
          \left(\frac{\log(u+1)}{L}\right)\right\}.
\end{equation}
Consequently,
\[
 |S_{j,k,\ell}(x,y;f)|
     \ll_{f,j,k,\ell}xL^{\kappa-1}\rho_\kappa(u).
\]
If \(C\ne0\), this bound has the correct order of magnitude, uniformly in
\eqref{eq:signed-range}. No sign hypothesis on the coefficients is imposed,
and no general nonvanishing assertion for odd \(\ell\) is needed.
\end{theorem}

\begin{theorem}[Zero prime mean]\label{thm:zero}
Suppose that \(\kappa=0\). There is \(c>0\), depending only on \(f,j,k,\ell\),
such that, uniformly for \(2\le y\le x\),
\begin{equation}\label{eq:zero-global}
 |S_{j,k,\ell}(x,y;f)|
   \ll_{f,j,k,\ell}
        x\{\log(2x)\}^{2D}\exp\{-c\sqrt{\log y}\}.
\end{equation}
There also exist constants \(A,c_1>0\), depending only on the same data, such
that throughout \eqref{eq:signed-range},
\begin{equation}\label{eq:zero-relative}
 |S_{j,k,\ell}(x,y;f)|
    \ll_{f,j,k,\ell}
        \Psi(x,y)\exp\{-c_1\sqrt{\log y}\}.
\end{equation}
\end{theorem}

An odd moment need not have zero prime mean. The following examples also give
nonvanishing constants in Theorem~\ref{thm:signed}.

\begin{corollary}[Two odd moments of sharp order]\label{cor:odd}
For every level-one form \(f\) as above,
\(\kappa_{2,2,1}=\kappa_{2,1,3}=1\) and both
\(C_{2,2,1}(f)\) and \(C_{2,1,3}(f)\) are positive. For a sufficiently large
\(A=A(f)>0\), uniformly in \eqref{eq:signed-range},
\begin{equation}\label{eq:odd-square}
 \sum_{\substack{n\le x\\P^+(n)\le y}}
       \lambda_{\Sym^2 f}(n^2)
   =C_{2,2,1}(f)x\rho(u)
       \left\{1+O_f\left(\frac{\log(u+1)}{L}\right)\right\},
\end{equation}
and
\begin{equation}\label{eq:odd-cubic}
 \sum_{\substack{n\le x\\P^+(n)\le y}}
       \lambda_{\Sym^2 f}(n)^3
   =C_{2,1,3}(f)x\rho(u)
       \left\{1+O_f\left(\frac{\log(u+1)}{L}\right)\right\}.
\end{equation}
Thus the upper order \(x\rho(u)\) is attained in both cases.
\end{corollary}

We next distinguish a second class of forms. Let \(f\) now be a fixed
normalized holomorphic CM newform of weight at least two, fixed level, and
trivial nebentypus. Only ordinary moments, \(j=k=1\), are considered in the
following theorem. Set
\[
 T_\ell(x,y;f)
   =\sum_{\substack{n\le x\\P^+(n)\le y}}\lambda_f(n)^\ell,
 \qquad
 \omega_{\CM}(y)
   =\frac{(\log y)^{3/5}}{\{\log\log(3y)\}^{1/5}}.
\]

\begin{theorem}[Ordinary moments of CM newforms]\label{thm:cm}
For \(\ell=2r\ge2\), put
\begin{equation}\label{eq:cm-index}
 \kappa_{\CM,2r}=\frac12\binom{2r}{r}.
\end{equation}
There is \(C^{\CM}_{f,2r}>0\) such that, for each fixed \(\varepsilon>0\),
uniformly in
\begin{equation}\label{eq:cm-even-range}
 \exp\{(\log\log x)^{5/3+\varepsilon}\}\le y\le x,
\end{equation}
one has
\begin{equation}\label{eq:cm-even}
 T_{2r}(x,y;f)
   =C^{\CM}_{f,2r}xL^{\kappa_{\CM,2r}-1}\rho_{\kappa_{\CM,2r}}(u)
      \left\{1+O_{f,r,\varepsilon}
          \left(\frac{\log(u+1)}{L}\right)\right\}.
\end{equation}
The constant is the prime-ordered product \eqref{eq:constant} with \(j=k=1\)
and \(\kappa=\kappa_{\CM,2r}\), including the actual bad-prime local factors.

For every odd \(\ell\ge1\), there is \(c>0\) such that
\begin{equation}\label{eq:cm-odd-global}
 |T_\ell(x,y;f)|
   \ll_{f,\ell}x\{\log(2x)\}^{2^{\ell+1}}
       \e^{-c\omega_{\CM}(y)}
       \qquad(2\le y\le x).
\end{equation}
There also exist \(A,c_1>0\), depending only on \(f,\ell\), such that, uniformly
in
\begin{equation}\label{eq:cm-odd-range}
 \exp\{A(\log x)^{5/8}(\log\log x)^{3/4}\}\le y\le x,
\end{equation}
one has
\begin{equation}\label{eq:cm-odd-relative}
 |T_\ell(x,y;f)|
    \ll_{f,\ell}\Psi(x,y)\e^{-c_1\omega_{\CM}(y)}.
\end{equation}
\end{theorem}

Absolute majorants also give first-moment estimates throughout the wider
nonnegative ranges.

\begin{corollary}[First moments in a wider range]\label{cor:first}
For \(m\ge1\), put
\begin{equation}\label{eq:mu}
 \mu_m=\frac2\pi\int_0^\pi|U_m(\cos\theta)|\sin^2\theta\,\mathrm d\theta
 =\frac{4(m+1)}{\pi m(m+2)}
       \cot\frac{\pi}{2(m+1)}<1.
\end{equation}
For a level-one \(f\), fixed \(m\ge1\), and \(0<\eta<1-\mu_m\), there is
\(A=A(f,m,\eta)>0\) such that, uniformly in \eqref{eq:even-range},
\begin{equation}\label{eq:first-noncm}
 \sum_{\substack{n\le x\\P^+(n)\le y}}
   \bigl\{|\lambda_{\Sym^m f}(n)|+|\lambda_f(n^m)|\bigr\}
    \ll_{f,m,\eta}\Psi(x,y)L^{\mu_m-1+\eta}.
\end{equation}
In particular, \(\mu_1=8/(3\pi)\).
For a fixed CM newform as in Theorem~\ref{thm:cm}, every
\(0<\eta<1-2/\pi\), and every \(\varepsilon>0\), uniformly in
\eqref{eq:cm-even-range},
\begin{equation}\label{eq:first-cm}
 \sum_{\substack{n\le x\\P^+(n)\le y}}|\lambda_f(n)|
       \ll_{f,\eta,\varepsilon}\Psi(x,y)L^{2/\pi-1+\eta}.
\end{equation}
The parameters \(m,\eta\) are fixed; no uniformity as \(\eta\downarrow0\) or
\(m\to\infty\) is asserted.
\end{corollary}

For each fixed positive \(\varepsilon\), the ranges
\eqref{eq:even-range}, \eqref{eq:signed-range}, and \eqref{eq:cm-odd-range}
contain, for sufficiently large \(x\), respectively
\[
 \begin{gathered}
 y\ge\exp\{(\log\log x)^{2+\varepsilon}\},\\
 y\ge\exp\{(\log x)^{2/3+\varepsilon}\},
 \qquad
 y\ge\exp\{(\log x)^{5/8+\varepsilon}\}.
 \end{gathered}
\]
The explicit logarithmic factors in the stated ranges follow by choosing
sufficiently large fixed constants in the error comparisons.

\subsection{Explicit indices and normalization checks}

Several consequences of the index computation are
\begin{align}
 \kappa_{j,1,2}=\kappa_{1,k,2}&=1,
 &\kappa_{j,1,4}&=j+1,\label{eq:indices-a}\\
 \kappa_{1,k,4}&=k+1,
 &\kappa_{j,2,2}&=\lfloor j/2\rfloor+1,\label{eq:indices-b}\\
 \kappa_{1,1,2r}&=\frac1{r+1}\binom{2r}{r}.
 &&\label{eq:indices-c}
\end{align}
Thus, for example, Theorem~\ref{thm:even} gives
\[
 \sum_{\substack{n\le x\\P^+(n)\le y}}
       \lambda_{\Sym^j f}(n^2)^2
   \sim C_{j,2,2}(f)xL^{\lfloor j/2\rfloor}
                    \rho_{\lfloor j/2\rfloor+1}(u)
\]
in its stated range. The second moment at \(n\) has index one, while the
second moment at \(n^2\) generally has a larger index.
The following exact values follow from \eqref{eq:index}.
\begin{table}[htbp]
\centering
\begin{tabular}{ccrrrr}
\toprule
\(j\)&\(k\)&\(\ell=1\)&\(\ell=2\)&\(\ell=3\)&\(\ell=4\)\\
\midrule
1&1&0&1&0&2\\
1&2&0&1&1&3\\
2&2&1&2&5&16\\
2&3&0&2&5&44\\
3&3&0&3&0&300\\
\bottomrule
\end{tabular}
\caption{Some values of the moment index \(\kappa_{j,k,\ell}\).}
\label{tab:indices}
\end{table}

For the ordinary second moment, the local identity
\[
 \sum_{\nu\ge0}\lambda_f(p^\nu)^2T^\nu
   =\frac{1+T}{(1-T)(1-\alpha_p^2T)(1-\alpha_p^{-2}T)}
\]
gives
\begin{equation}\label{eq:second-moment-check}
 \sum_{n\ge1}\frac{\lambda_f(n)^2}{n^s}
   =\frac{\zeta(s)L(s,\Sym^2f)}{\zeta(2s)},\qquad
 C_{1,1,2}(f)=\frac{L(1,\Sym^2f)}{\zeta(2)}.
\end{equation}
For example, the local identity follows by squaring
\(\lambda_f(p^\nu)=U_\nu(\cos\theta_p)\) and summing the resulting geometric
series; continuity covers \(\alpha_p=\pm1\).
For \(\kappa\ge1\), at \(y=x\) one has \(u=1\) and
\(\rho_\kappa(1)=1/\Gamma(\kappa)\), so the main term is
\[
 \frac{C_{j,k,\ell}(f)}{\Gamma(\kappa)}
       x(\log x)^{\kappa-1}.
\]
This agrees with the leading residue coefficient in
Lemma~\ref{lem:unrestricted}.

\section{The local representation and the moment index}
\label{sec:representation}

For a finite-dimensional complex representation $E$, write $\Sym^r E$
for its $r$th symmetric tensor power and $E^{\otimes r}$ for its $r$th
tensor power. The notation $E^{\oplus b}$ means the direct sum of $b$
copies of $E$, with $E^{\oplus0}=0$. Put $V_m=\Sym^m\C^2$, where
$\C^2$ carries the standard action of $\SL_2(\C)$. The character of a
representation is the trace of its action. On the diagonal torus the
character of $V_m$ is
\[
  \chi_m(z)=z^m+z^{m-2}+\cdots+z^{-m},
  \qquad \chi_m(\e^{i\theta})=U_m(\cos\theta).
\]
The relationship between coefficients at perfect powers and symmetric-power
representations is already recorded in \cite[Section~2.1]{YY}. For $k=1$,
the prime decompositions below agree with those in \cite[Lemma~2.2]{Venk}
and \cite[Sections~2.3--2.4]{YY}; in particular, our
$\kappa_{j,1,\ell}$ is the multiplicity denoted by $K_{0,j,\ell}$ in the
latter paper. We use the same representation-theoretic principle for the
plethysm $\Sym^k(V_j)$, retaining its effect on the friable moment problem.
For background on these plethysms and Hermite reciprocity, see \cite{PW}.

\begin{lemma}\label{lem:representation}
Let $W_{j,k}=\Sym^k(V_j)$. There are nonnegative integers $b_m$, indexed by
$0\le m\le jk\ell$, such that
\begin{equation}\label{eq:decomposition}
  W_{j,k}^{\otimes\ell}
  \simeq \bigoplus_{m=0}^{jk\ell}V_m^{\oplus b_m}.
\end{equation}
At every prime $p$,
\begin{equation}\label{eq:prime-decomposition}
  \lambda_{\Sym^j f}(p^k)^\ell
  =\sum_{m=0}^{jk\ell}b_m\lambda_{\Sym^m f}(p),
  \qquad \lambda_{\Sym^0 f}(p)=1.
\end{equation}
Moreover, $b_0=\kappa$ and
\begin{equation}\label{eq:degree-sum}
  \sum_{m=0}^{jk\ell}b_m(m+1)=D.
\end{equation}
We extend the multiplicities by $b_m=0$ for every integer $m>jk\ell$.
\end{lemma}

\begin{proof}
The coefficient of $T^k$ in \eqref{eq:local} is the complete homogeneous
symmetric polynomial of degree $k$ in
$\alpha_p^j,\alpha_p^{j-2},\ldots,\alpha_p^{-j}$. It is the character of
$W_{j,k}$ at $\operatorname{diag}(\alpha_p,\alpha_p^{-1})$, namely
$Q_{j,k}(\alpha_p)$. Complete reducibility for finite-dimensional
$\SL_2(\C)$-representations gives \eqref{eq:decomposition}, and evaluation
of its character gives \eqref{eq:prime-decomposition}.

Write $Q_{j,k}(z)^\ell=\sum_r c_rz^r$, with $c_r=0$ outside its finite
support. Each $V_m$ has weights $m,m-2,\ldots,-m$, each with multiplicity
one. Consequently,
\[
  c_m=b_m+b_{m+2}+b_{m+4}+\cdots,
  \qquad b_m=c_m-c_{m+2}\quad(m\ge0).
\]
Taking $m=0$ gives precisely the definition of $\kappa$ in
\eqref{eq:index}. Finally,
$\dim W_{j,k}=\binom{j+k}{j}$, so taking dimensions in
\eqref{eq:decomposition} gives \eqref{eq:degree-sum}.
\end{proof}

Let $\mu_{\ST}$ be the Sato--Tate probability measure on $[0,\pi]$,
with $d\mu_{\ST}(\theta)=(2/\pi)\sin^2\theta\,d\theta$. Equivalently,
\begin{equation}\label{eq:index-integral}
  \kappa=\int_0^\pi Q_{j,k}(\e^{i\theta})^\ell\,d\mu_{\ST}(\theta).
\end{equation}
Indeed, $\chi_0=1$, and the character formula gives
\[
 \begin{aligned}
 \int_0^\pi\chi_m(\e^{i\theta})\chi_n(\e^{i\theta})
       \,d\mu_{\ST}(\theta)
 &=\frac2\pi\int_0^\pi\sin((m+1)\theta)\sin((n+1)\theta)\,d\theta\\
 &=\begin{cases}1,&m=n,\\0,&m\ne n.\end{cases}
 \end{aligned}
\]
Applying this with $n=0$ to \eqref{eq:decomposition} proves
\eqref{eq:index-integral}. This is an identity between characters; the effective
prime-number theorem used later is a separate analytic input.

\begin{proposition}\label{prop:parity}
If $jk\ell$ is odd, then $\kappa=0$. For $\ell\ge2$, one has
$\kappa>0$ if and only if $jk\ell$ is even. For $\ell=1$, evenness of
$jk$ does not imply positivity.
\end{proposition}

\begin{proof}
All the weights of $W_{j,k}^{\otimes\ell}$ have the parity of $jk\ell$.
If this number is odd, the trivial representation cannot occur.
Every $V_m$ is self-dual. Since $W_{j,k}$ is a nonzero direct sum of such
representations, $W_{j,k}^{\otimes2}$ has a nonzero invariant. Pairing
tensor factors therefore proves positivity for every even $\ell$.

Suppose next that $\ell\ge3$ is odd and $jk$ is even. The highest weight
$jk$ occurs in $W_{j,k}$, so $V_{jk}$ is a constituent. The Clebsch--Gordan
identity
\begin{equation}\label{eq:clebsch-gordan}
  V_m\otimes V_n
  \simeq\bigoplus_{r=0}^{\min(m,n)}V_{m+n-2r}
\end{equation}
shows that $V_{jk}$ occurs in $V_{jk}\otimes V_{jk}$: take
$r=jk/2$. Self-duality then implies that $V_{jk}^{\otimes3}$ contains
an invariant. Thus $W_{j,k}^{\otimes3}$ has an invariant, and pairing the
remaining $\ell-3$ factors proves the assertion. Finally,
$W_{1,k}=V_k$ has no invariant for $k\ge1$, including even $k$.
\end{proof}

\subsection{Coefficient formulas and examples}

Let
\[
  \mathcal G_{j,k}(q)
  :=\begin{bmatrix}j+k\\k\end{bmatrix}_{q}
  =\sum_{\substack{a_0+\cdots+a_j=k\\a_r\ge0}}
       q^{\sum_{r=0}^j r a_r}
\]
be the Gaussian binomial polynomial. Since
$Q_{j,k}(z)=z^{jk}\mathcal G_{j,k}(z^{-2})$, the index can also be written
as
\begin{equation}\label{eq:gaussian-index}
  \kappa
  =[q^{jk\ell/2}]\mathcal G_{j,k}(q)^\ell
   -[q^{jk\ell/2-1}]\mathcal G_{j,k}(q)^\ell
  \qquad(2\mid jk\ell).
\end{equation}
This gives a finite procedure using only integer polynomial arithmetic.
The polynomial $\mathcal G_{j,k}$ counts partitions into at most $k$
parts, each at most $j$, with the exponent recording the size of the
partition. Conjugating their diagrams gives
$\mathcal G_{j,k}(q)=\mathcal G_{k,j}(q)$ and hence
\begin{equation}\label{eq:index-reciprocity}
  \kappa_{j,k,\ell}=\kappa_{k,j,\ell}.
\end{equation}
The corresponding symmetry of representations is Hermite reciprocity;
see \cite{PW}. The symmetry of the indices does not in general extend to
the arithmetic constants. An exact comparison for $(j,k,\ell)=(2,1,1)$
and $(1,2,1)$ is proved after Lemma~\ref{lem:euler}; see
\eqref{eq:constant-asymmetry}.

We now verify \eqref{eq:indices-a}--\eqref{eq:indices-c}. The second-
and fourth-moment identities with $j=1$ or $k=1$ follow from
$W_{j,1}=V_j$, $W_{1,k}=V_k$ and character orthogonality.
In particular, \eqref{eq:clebsch-gordan} decomposes
$V_j\otimes V_j$ into $j+1$ distinct irreducibles, so the squared norm of
its character, and hence the fourth moment of $\chi_j$, is $j+1$.
For $\kappa_{j,2,2}$, use the multiplicity-free decomposition
\begin{equation}\label{eq:symmetric-square-decomposition}
  \Sym^2(V_j)\simeq
  \bigoplus_{r=0}^{\lfloor j/2\rfloor}V_{2j-4r}.
\end{equation}
To verify this decomposition, the character of $\Sym^2(V_j)$ is
$\{\chi_j(z)^2+\chi_j(z^2)\}/2$. For $0\le r\le j$, the coefficient
of $z^{2j-2r}$ is $\lfloor r/2\rfloor+1$: the first summand contributes
$(r+1)/2$, and the second contributes $1/2$ exactly when $r$ is even.
Taking successive differences of these weight multiplicities, as in
Lemma~\ref{lem:representation}, gives multiplicity one for
$V_{2j-4r}$, $0\le r\le\lfloor j/2\rfloor$, and zero for the remaining
irreducibles. This proves \eqref{eq:symmetric-square-decomposition}.
The squared norm of the character is therefore $\lfloor j/2\rfloor+1$.
Finally, the coefficients of $z^0$ and $z^2$ in
$(z+z^{-1})^{2r}$ are $\binom{2r}{r}$ and $\binom{2r}{r-1}$,
respectively. Their difference is the Catalan number in
\eqref{eq:indices-c}.

Table~\ref{tab:indices} is obtained from \eqref{eq:gaussian-index} or by
successive applications of \eqref{eq:clebsch-gordan}. For the last three
rows the starting decompositions are
\[
  W_{2,2}\simeq V_4\oplus V_0,\qquad
  W_{2,3}\simeq V_6\oplus V_2,\qquad
  W_{3,3}\simeq V_9\oplus V_5\oplus V_3.
\]
Their characters can be checked directly from \eqref{eq:Q}; multiplying
these Laurent polynomials and taking the difference of the coefficients
of $z^0$ and $z^2$ gives the displayed entries. No numerical
Sato--Tate integration is involved.

\section{Arithmetic and analytic input}
\label{sec:arithmetic}

Throughout this section set
\[
  a(n)=\lambda_{\Sym^j f}(n^k)^\ell.
\]
The integers $j,k,\ell$ and the level-one form $f$ are fixed.

\begin{lemma}[A divisor majorant]\label{lem:majorant}
The function $a$ is real and multiplicative, with $a(1)=1$, and
\begin{equation}\label{eq:divisor-majorant}
  |a(n)|\le d_D(n),\qquad D=\binom{j+k}{j}^{\ell}.
\end{equation}
In particular,
\[
  |a(p^\nu)|\ll_{j,k,\ell}(\nu+1)^{j\ell},
  \qquad
  \sum_{n\le t}\frac{|a(n)|}{n}\ll_D\{\log(2t)\}^{D}.
\]
\end{lemma}

\begin{proof}
The local parameters occur in conjugate pairs, so all the coefficients
are real. Multiplicativity follows by applying Hecke multiplicativity to
$m^k,n^k$ when $(m,n)=1$. Equation~\eqref{eq:local} and
$|\alpha_p|=1$ give
\[
  |a(p^\nu)|\le\binom{k\nu+j}{j}^{\ell}.
\]
Put $E=\Sym^k\C^{j+1}$, so that $\dim E=\binom{j+k}{j}$. Multiplication
of homogeneous polynomials in each tensor position defines a map
\begin{equation}\label{eq:majorant-surjection}
  \Sym^\nu\!\left(E^{\otimes\ell}\right)
  \longrightarrow
  \left(\Sym^{k\nu}\C^{j+1}\right)^{\otimes\ell}.
\end{equation}
It is well defined because multiplication is invariant under simultaneous
permutation of the $\nu$ inputs. It is surjective: a monomial of degree
$k\nu$ in each of the $\ell$ tensor positions can be split into $\nu$
monomials of degree $k$ in that position, and grouping the resulting
factors produces a preimage. Comparing dimensions in
\eqref{eq:majorant-surjection} gives
\[
  \binom{k\nu+j}{j}^{\ell}
  \le \binom{D+\nu-1}{\nu}=d_D(p^\nu).
\]
Multiplicativity proves \eqref{eq:divisor-majorant}, and the polynomial
bound at prime powers follows from the preceding binomial estimate.
Finally, positivity of the divisor coefficients and Mertens' estimate
give
\[
  \sum_{n\le t}\frac{|a(n)|}{n}
  \le\sum_{n\le t}\frac{d_D(n)}{n}
  \le\prod_{p\le t}(1-p^{-1})^{-D}
  \ll_D\{\log(2t)\}^{D}.
\]
\end{proof}

\begin{lemma}[Prime mean]\label{lem:prime}
There is $c=c(f,j,k,\ell)>0$ such that
\begin{equation}\label{eq:prime-mean}
  \sum_{p\le t}a(p)\log p
  =\kappa t+O_{f,j,k,\ell}\!\left(t\e^{-c\sqrt{\log t}}\right)
  \qquad(t\ge3).
\end{equation}
\end{lemma}

\begin{proof}
For each fixed $m\ge1$, \cite[Proposition~2.1]{Thorner} gives
\[
  \sum_{p\le t}U_m(\cos\theta_p)\log p
  \ll_{f,m}t\e^{-c_{f,m}\sqrt{\log t}}
\]
with $c_{f,m}>0$. The degree restriction in that proposition holds for all
sufficiently large $t$ when $m$ is fixed, and the remaining bounded range
is absorbed into the implied constant. All the symmetric powers required
here are cuspidal automorphic by \cite{NT}. Apply
\eqref{eq:prime-decomposition}, use
$\lambda_{\Sym^m f}(p)=U_m(\cos\theta_p)$, and use the classical
prime-number theorem for $m=0$. Only finitely many $m$ occur, and their
trivial multiplicity is $b_0=\kappa$.
\end{proof}

\begin{lemma}[Euler factorization and the constant]\label{lem:euler}
For $\Re s>1$,
\begin{equation}\label{eq:factorization}
\begin{split}
  F(s):=\sum_{n\ge1}\frac{a(n)}{n^s}&=\zeta(s)^\kappa G(s),\\
  G(s)&=H(s)\prod_{m=1}^{jk\ell}L(s,\Sym^m f)^{b_m}.
\end{split}
\end{equation}
The Euler product for $H$ converges absolutely and locally uniformly in
$\Re s>1/2$, and $H$ is bounded on each half-plane
$\Re s\ge1/2+\eta$, $\eta>0$. The function $G$ is holomorphic in
$\Re s>1/2$. The prime-ordered limit in \eqref{eq:constant} exists and
equals $G(1)$. We use $C_{j,k,\ell}(f)$ for this value also when
$\kappa=0$. If every local sum
$\sum_{\nu\ge0}a(p^\nu)p^{-\nu}$ is positive, then $G(1)>0$;
in particular this holds when $\ell$ is even.
\end{lemma}

\begin{proof}
Let
\[
  A_p(s)=\sum_{\nu\ge0}a(p^\nu)p^{-\nu s},\qquad
  H_p(s)=A_p(s)\prod_{m=0}^{jk\ell}L_p(s,\Sym^m f)^{-b_m},
\]
where, for every $m\ge0$,
\[
 L_p(s,\Sym^m f)
  =\prod_{r=0}^{m}(1-\alpha_p^{m-2r}p^{-s})^{-1}.
\]
In particular, $L_p(s,\Sym^0 f)=(1-p^{-s})^{-1}$.
The coefficient of $p^{-s}$
in $H_p(s)$ vanishes by \eqref{eq:prime-decomposition}. The polynomial
bound from Lemma~\ref{lem:majorant} therefore gives
\begin{equation}\label{eq:H-local-bound}
  H_p(s)=1+O_{j,k,\ell,\eta}(p^{-2\Re s})
  \qquad(\Re s\ge1/2+\eta).
\end{equation}
Every local series $A_p$ is holomorphic for $\Re s>0$, and the other
local factors in $H_p$ are polynomials in $p^{-s}$. Thus
\eqref{eq:H-local-bound} proves the absolute and locally uniform
convergence of $H(s)=\prod_pH_p(s)$, as well as the stated boundedness.
Equation~\eqref{eq:factorization} follows first in $\Re s>1$ by absolute
convergence. The nontrivial standard symmetric-power $L$-functions are
entire, by \cite{NT} and the analytic theory of cuspidal standard
$L$-functions, so $G$ is holomorphic in $\Re s>1/2$.

We next prove convergence of the boundary Euler product. Put
$b(p)=a(p)-\kappa$ and
$B(t)=\sum_{p\le t}b(p)\log p$. Lemma~\ref{lem:prime} and the classical
prime-number theorem imply, for some $c_1>0$,
\begin{equation}\label{eq:boundary-prime-sum}
  B(t)=\sum_{p\le t}a(p)\log p
       -\kappa\sum_{p\le t}\log p
  \ll t\e^{-c_1\sqrt{\log t}}\qquad(t\ge3).
\end{equation}
For $Z>z\ge3$ and real $1\le\sigma\le2$, partial summation gives
\begin{align*}
  \sum_{z<p\le Z}\frac{b(p)}{p^\sigma}
  ={}&\frac{B(Z)}{Z^\sigma\log Z}
       -\frac{B(z)}{z^\sigma\log z}\\
     &+\int_z^Z\frac{B(t)}{t^{\sigma+1}}
       \left\{\frac{\sigma}{\log t}
                    +\frac{1}{(\log t)^2}\right\}\,dt.
\end{align*}
The first term tends to zero uniformly in $1\le\sigma\le2$ as
$Z\to\infty$. The integral converges absolutely and uniformly there,
since its absolute value is bounded by a constant multiple of
\[
  \int_z^\infty\frac{\e^{-c_1\sqrt{\log t}}}{t\log t}\,dt
  =2\int_{\sqrt{\log z}}^\infty\frac{\e^{-c_1v}}{v}\,dv
  \ll\frac{\e^{-c_1\sqrt{\log z}}}{\sqrt{\log z}}.
\]
The boundary term at $z$ is smaller. Decreasing the positive constant
in the exponent if necessary, we obtain
\begin{equation}\label{eq:boundary-linear-tail}
  \sum_{p>z}\frac{b(p)}{p^\sigma}
  \ll\e^{-c_2\sqrt{\log z}}
  \qquad(1\le\sigma\le2)
\end{equation}
uniformly. In particular, $\sum_p b(p)/p$ converges in increasing prime
order.

For real $\sigma\ge1$, write
$E_p(\sigma)=(1-p^{-\sigma})^\kappa A_p(\sigma)$. The local polynomial
bound gives, uniformly for $\sigma\ge1$,
\[
  \sum_{\nu\ge2}|a(p^\nu)|p^{-\nu\sigma}\ll p^{-2\sigma}.
\]
Consequently,
\begin{equation}\label{eq:boundary-local-expansion}
  E_p(\sigma)=1+b(p)p^{-\sigma}+O(p^{-2\sigma}).
\end{equation}
Since $b(p)$ is bounded, there is a fixed $p_0\ge3$ such that
$|E_p(\sigma)-1|\le1/2$ for $p>p_0$ and $\sigma\ge1$. These real
factors are positive, and
\begin{equation}\label{eq:boundary-log-expansion}
  \log E_p(\sigma)=b(p)p^{-\sigma}+O(p^{-2\sigma}).
\end{equation}
The error series in \eqref{eq:boundary-log-expansion} converges absolutely
and uniformly for $1\le\sigma\le2$. Together with
\eqref{eq:boundary-linear-tail}, this proves uniform convergence of
$\sum_{p>p_0}\log E_p(\sigma)$ on that interval. The prime-ordered
product in \eqref{eq:constant} therefore exists and equals
\begin{equation}\label{eq:constant-finite-tail}
  \left(\prod_{p\le p_0}E_p(1)\right)
  \exp\left\{\sum_{p>p_0}\log E_p(1)\right\}.
\end{equation}
All the finitely many small-prime factors have been kept outside the
logarithm. In a signed moment they may be negative or zero; a zero factor
makes every sufficiently long partial product zero. The exponential
tail in \eqref{eq:constant-finite-tail} is always strictly positive.

For $\sigma>1$, absolute convergence gives
\[
  G(\sigma)=\prod_pE_p(\sigma)
  =\left(\prod_{p\le p_0}E_p(\sigma)\right)
   \exp\left\{\sum_{p>p_0}\log E_p(\sigma)\right\}.
\]
The finite product is continuous at $\sigma=1$, and the tail is continuous
there by its uniform convergence. Letting $\sigma\downarrow1$ and using
the holomorphy of $G$ identifies \eqref{eq:constant-finite-tail} with
$G(1)$. Thus boundary convergence has been proved directly, rather than
inferred from analytic continuation. If $A_p(1)>0$ for every prime, then
all the finite factors in \eqref{eq:constant-finite-tail} are positive,
so $G(1)>0$. For even $\ell$, all $a(p^\nu)$ are nonnegative and
$A_p(1)\ge1$, which proves the last assertion.
\end{proof}

\paragraph{An explicit asymmetry of the constants.}
Both $(j,k,\ell)=(2,1,1)$ and $(1,2,1)$ have $\kappa=0$, since their
prime character is $\chi_2$. The local generating functions are
\[
 \begin{split}
 \mathcal A_p(T)
  &:=\sum_{\nu\ge0}\lambda_{\Sym^2 f}(p^\nu)T^\nu
    =\frac1{(1-T)(1-\alpha_p^2T)(1-\alpha_p^{-2}T)},\\
 \mathcal B_p(T)
  &:=\sum_{\nu\ge0}\lambda_f(p^{2\nu})T^\nu
    =\frac{1+T}{(1-\alpha_p^2T)(1-\alpha_p^{-2}T)}
    =(1-T^2)\mathcal A_p(T).
 \end{split}
\]
The second identity follows by extracting the even coefficients in
$(1-\alpha_p z)^{-1}(1-\alpha_p^{-1}z)^{-1}$: average its values at
$z$ and $-z$, and put $T=z^2$. Thus, by absolute convergence for
$\Re s>1$,
\[
 \sum_{n\ge1}\frac{\lambda_f(n^2)}{n^s}
   =\frac{L(s,\Sym^2 f)}{\zeta(2s)}.
\]
Moreover,
$\mathcal A_p(p^{-1})=(1-p^{-1})^{-1}|1-\alpha_p^2p^{-1}|^{-2}>0$,
and $\mathcal B_p(p^{-1})=(1-p^{-2})\mathcal A_p(p^{-1})>0$.
Lemma~\ref{lem:euler} therefore gives
\begin{equation}\label{eq:constant-asymmetry}
  C_{2,1,1}(f)=L(1,\Sym^2 f)
     =\zeta(2)C_{1,2,1}(f)>C_{1,2,1}(f)>0.
\end{equation}
This proves the failure of symmetry for the complete constants, rather
than only a difference between their individual local factors.

\begin{lemma}[Unrestricted summation and an analytic class]
\label{lem:unrestricted}
There is $\delta_0>0$ such that
\begin{equation}\label{eq:unrestricted}
  A(t):=\sum_{n\le t}a(n)
  =\Res_{s=1}\frac{F(s)t^s}{s}+O(t^{1-\delta_0}).
\end{equation}
If $\kappa=0$, the residue is absent. If $\kappa\ge1$, it is
$tP(\log t)$, where $\deg P\le\kappa-1$ and the coefficient of
$(\log t)^{\kappa-1}$ is $G(1)/\Gamma(\kappa)$. Furthermore, for some
$\delta>0$,
\begin{equation}\label{eq:analytic-class}
  G\text{ is holomorphic on }\Re s>1-\delta,
  \qquad
  G(s)\ll(1+|\Im s|)^{1-\delta}\quad(\Re s>1-\delta).
\end{equation}
\end{lemma}

\begin{proof}
The degree identity \eqref{eq:degree-sum}, the boundedness of $H$, and
convexity for the fixed standard $L$-functions give
\begin{equation}\label{eq:F-convexity}
  F(\sigma+i\tau)
  \ll_{\eta,\varepsilon}(|\tau|+2)^{D(1-\sigma)/2+\varepsilon}
\end{equation}
for $1/2+\eta\le\sigma\le1$ and $|\tau|\ge2$; see
\cite[Chapter~5]{IK}. In $1\le\sigma\le1+1/\log t_0$ the corresponding
bound is $F(\sigma+i\tau)\ll_\varepsilon(|\tau|+2)^\varepsilon$.

First take a half-integer endpoint $t_0\ge3$, and put
$c_0=1+1/\log t_0$ and $T=t_0^{2/(D+4)}$. Truncated Perron inversion,
using $|a(n)|\le d_D(n)\ll_\varepsilon n^\varepsilon$ and the distance
of $t_0$ from the integers, gives
\[
  A(t_0)=\frac{1}{2\pi i}\int_{c_0-iT}^{c_0+iT}
                F(s)\frac{t_0^s}{s}\,ds
           +O_\varepsilon\left(\frac{t_0^{1+\varepsilon}}{T}\right).
\]
Any logarithmic factors in the usual truncation estimate are absorbed by
decreasing the preliminary coefficient exponent and then renaming
$\varepsilon$. Move the contour to $\sigma_0=1/2+\eta$, where
$0<\eta<1/4$. The only possible pole crossed is at $s=1$.
The part of the new vertical integral with bounded imaginary part is
$O(t_0^{\sigma_0})$. For the remaining part,
\eqref{eq:F-convexity} gives the bound
\[
  O_{\eta,\varepsilon}\left(
       t_0^{1/2+\eta}T^{D(1/2-\eta)/2+\varepsilon}\right).
\]
On either horizontal side, in the portion $\sigma_0\le\sigma\le1$,
the quantity $t_0^\sigma T^{D(1-\sigma)/2}$ increases with $\sigma$,
because
\[
  \log t_0-\frac{D}{2}\log T=\frac{4}{D+4}\log t_0>0.
\]
In the remaining portion $1\le\sigma\le c_0$, use
$t_0^\sigma\le\e t_0$ and the bound $F(\sigma+iT)\ll T^\varepsilon$.
The horizontal contributions are therefore
$O_\varepsilon(t_0T^{-1+\varepsilon})$.

The exponent in the vertical bound is
\[
  \frac{D+2+4\eta}{D+4}+\frac{2\varepsilon}{D+4},
\]
and the Perron error has exponent $1-2/(D+4)+\varepsilon$.
Choosing $\eta$ and then $\varepsilon$ sufficiently small makes both
strictly less than one. The residue theorem proves
\eqref{eq:unrestricted} at half-integer endpoints for some fixed
$\delta_0>0$. An arbitrary endpoint differs from a suitable half-integer
by at most one; its summatory function changes by at most
$O_\varepsilon(t^\varepsilon)$ and the residue polynomial changes by
a fixed power of $\log(2t)$. Reducing $\delta_0$ if necessary proves
\eqref{eq:unrestricted} for all $t\ge3$.

If $\kappa=0$, $F$ is holomorphic at one. If $\kappa\ge1$, the Laurent
expansion of $\zeta(s)^\kappa G(s)t^s/s$ shows that its residue is
$tP(\log t)$ with degree at most $\kappa-1$; the coefficient of the
highest possible degree is $G(1)/(\kappa-1)!=G(1)/\Gamma(\kappa)$.

Finally, the total degree of the nontrivial $L$-factors of $G$ is
$D-\kappa$. Choose
\[
  0<\delta<\min\left\{\frac14,\frac{1}{2D+2}\right\}.
\]
Convexity in $1-\delta\le\sigma\le1$, with a sufficiently small loss,
gives
$G(\sigma+i\tau)\ll(|\tau|+2)^{(D-\kappa)\delta/2+\varepsilon}$,
whose exponent is less than $1-\delta$. For $\sigma\ge1$, the same
boundary bound and absolute convergence to its right give the required
estimate. Bounded imaginary parts are covered by holomorphy and
compactness on each fixed strip, and by absolute convergence farther
to the right. Lemma~\ref{lem:euler} supplies holomorphy on this
half-plane, proving \eqref{eq:analytic-class}.
\end{proof}

\section{An exponential convolution and prime replacement}
\label{sec:replacement}

We first isolate an elementary passage from an effective prime mean to a
friable sum.  No sign condition is imposed on the coefficients in this
section.  For a real number $\kappa\geq 0$, let $z_\kappa$ be the function
specified by
\begin{equation}
\label{eq:z-kappa}
\begin{split}
 z_\kappa(v)&=0 \quad (v<0),\\
 z_\kappa(v)&=1 \quad (0\leq v\leq 1),\\
 v z_\kappa'(v)&=-\kappa z_\kappa(v-1) \quad (v>1).
\end{split}
\end{equation}
The differential--difference equation determines the function successively
on the intervals $(r,r+1]$.  Given
\[
 A(t)=\sum_{n\leq t}a(n),
\]
with $A(t)=0$ for $0<t<1$, define
\begin{equation}
\label{eq:abstract-main}
 \Lambda_a(x,y)
   =x\int_{\R}z_\kappa(u-v)\,
       d\!\left\{\frac{A(y^v)}{y^v}\right\}.
\end{equation}
The dependence of $\Lambda_a$ on the fixed parameter $\kappa$ is suppressed
in the notation.  The functions in this Stieltjes integral are taken to be
right-continuous.
In particular, the measure retains its atom at $v=0$.  The full summatory
function $A$, rather than an asymptotic approximation to it, occurs in
\eqref{eq:abstract-main}.

\begin{lemma}[Prime replacement]
\label{lem:replacement}
Let $D\geq 1$ be an integer, and let $a$ be a real multiplicative function
with $a(1)=1$ and $|a(n)|\leq d_D(n)$.  Suppose that $0\leq\kappa\leq D$
and that, for some $c>0$,
\begin{equation}
\label{eq:replacement-prime-hypothesis}
 \sum_{p\leq t}a(p)\log p
   =\kappa t+O\bigl(t\e^{-c\omega(t)}\bigr)
       \qquad(t\geq 3),
\end{equation}
where either $\omega(t)=\sqrt{\log t}$ or
$\omega(t)=\omega_{\CM}(t)$.  Then there is $c_1>0$ such that, uniformly
for $2\leq y\leq x$,
\begin{equation}
\label{eq:replacement}
 \sum_{\substack{n\leq x\\P^+(n)\leq y}}a(n)
   =\Lambda_a(x,y)
     +O\!\left(x\{\log(2x)\}^{2D}\e^{-c_1\omega(y)}\right).
\end{equation}
The constants may depend on $D$ and on the constants in
\eqref{eq:replacement-prime-hypothesis}.  If $\kappa=0$, then
$\Lambda_a(x,y)=A(x)$.
\end{lemma}

\begin{proof}
Define a multiplicative function $e_y$ by $e_y(1)=1$ and, for $\nu\geq 1$,
\[
 e_y(p^\nu)=
 \begin{cases}
  (-a(p))^\nu/\nu!,&p>y,\\
  0,&p\leq y,
 \end{cases}
\]
and put $c_y=a*e_y$.  Since $|a(p)|\leq D$, comparison of the
nonnegative coefficients of $\exp(DT)$ and $(1-T)^{-D}$ gives
$|e_y|\leq d_D$.  Consequently
\[
 |c_y|\leq d_D*d_D=d_{2D}.
\]
We have $c_y(n)=a(n)$ whenever $P^+(n)\leq y$, whereas $c_y(p)=0$ for
$p>y$.  Thus a nonzero value $c_y(n)$ with $P^+(n)>y$ requires a prime
$p>y$ for which $p^2\mid n$.  The elementary inequalities
$d_{2D}(p^2m)\leq d_{2D}(p^2)d_{2D}(m)$ and
$\sum_{m\leq z}d_{2D}(m)\ll_D z\{\log(2z)\}^{2D-1}$ therefore imply
\begin{equation}
\label{eq:squareful-error}
\begin{split}
 \left|\sum_{n\leq x}c_y(n)
       -\sum_{\substack{n\leq x\\P^+(n)\leq y}}a(n)\right|
 &\leq \sum_{p>y}\ \sum_{m\leq x/p^2}d_{2D}(p^2m)\\
 &\ll_D \frac{x\{\log(2x)\}^{2D-1}}{y}.
\end{split}
\end{equation}

Expansion of the exponential Euler factors gives the exact identity
\begin{equation}
\label{eq:exponential-identity}
 \sum_{n\leq x}c_y(n)
   =\sum_{r\geq 0}\frac{(-1)^r}{r!}
       \sum_{\substack{mp_1\cdots p_r\leq x\\p_1,\ldots,p_r>y}}
        a(m)a(p_1)\cdots a(p_r).
\end{equation}
The primes are ordered and may repeat.  The factorial denominators in
$e_y(p^\nu)$ are exactly those obtained by counting the orderings of a
multiset of primes.  There is no coprimality restriction between $m$ and
the $p_i$.  For fixed $x,y$, only finitely many $r$ contribute.

Partial summation in \eqref{eq:replacement-prime-hypothesis} gives, for
$Z\geq y\geq 2$,
\begin{equation}
\label{eq:unweighted-prime-replacement}
 \left|\sum_{y<p\leq Z}a(p)-\kappa\int_y^Z\frac{dt}{\log t}\right|
   \ll Z\e^{-c_2\omega(y)}
\end{equation}
with a suitable $c_2>0$.  Indeed, the endpoint errors and the integral
of the weighted error are bounded by this expression; the finitely
bounded interval below $3$ can be absorbed in the implied constant.
Both choices of $\omega$ are increasing for all sufficiently large
arguments, and decreasing $c_2$ deals with the remaining bounded range.

Replace the $r$ prime measures in \eqref{eq:exponential-identity}, one
at a time, by $\kappa\,dt/\log t$.  After fixing $m$ and the other
$r-1$ variables, the variable being replaced ranges over $(y,Z]$, where
\[
 Z=\frac{x}{mt_1\cdots t_{r-1}}.
\]
If $Z<y$, both its discrete and its continuous contribution vanish.
Otherwise \eqref{eq:unweighted-prime-replacement} applies directly.
This procedure needs no estimate for the variation of $A$.

The harmonic total masses of each of the remaining discrete or
continuous measures are at most
\[
 M=D\log\log(2x)+O_D(1),
\]
because
\[
 \sum_{p\leq x}\frac{|a(p)|}{p}
      \leq D\log\log(2x)+O_D(1),
 \qquad
 \kappa\int_y^x\frac{dt}{t\log t}
      \leq D\log\log(2x)+O_D(1).
\]
Increase the constant so that $M\geq 1$.  Taking absolute values of all
other measures, and then enlarging their ranges, bounds the error in
the $r$th term before division by $r!$ by
\[
 C r x\e^{-c_2\omega(y)}
    \left(\sum_{m\leq x}\frac{|a(m)|}{m}\right)M^{r-1}.
\]
The divisor majorant gives
$\sum_{m\leq x}|a(m)|/m\ll_D\{\log(2x)\}^D$.  Since
\[
 \sum_{r\geq 1}\frac{rM^{r-1}}{r!}=\e^M
       \ll_D\{\log(2x)\}^D,
\]
the total replacement error is
\begin{equation}
\label{eq:total-replacement-error}
 \ll x\{\log(2x)\}^{2D}\e^{-c_2\omega(y)}.
\end{equation}
The resulting continuous expression is
\begin{equation}
\label{eq:continuous-replacement}
 M_a(x,y)=\sum_{r\geq 0}\frac{(-\kappa)^r}{r!}
   \int_{\substack{t_1,\ldots,t_r>y\\t_1\cdots t_r\leq x}}
     A\!\left(\frac{x}{t_1\cdots t_r}\right)
          \prod_{i=1}^r\frac{dt_i}{\log t_i},
\end{equation}
where the term $r=0$ is $A(x)$.  This sum is also finite for fixed $x,y$.

It remains to identify \eqref{eq:continuous-replacement} exactly with
\eqref{eq:abstract-main}.  Put
\[
 h_y(v)=\frac{A(y^v)}{y^v}\quad(v\geq 0),\qquad
 h_y(v)=0\quad(v<0),\qquad
 h(v)=
 \begin{cases}
  1/v,&v\geq 1,\\
  0,&v<1.
 \end{cases}
\]
Let $\mathfrak d_0$ denote the Dirac unit mass at zero, and let $h^{*r}$
be the $r$-fold additive convolution of $h$.  The substitutions
$t_i=y^{v_i}$ give
\begin{equation}
\label{eq:kernel-convolution}
 \frac{M_a(x,y)}{x}=(K_\kappa*h_y)(u),\qquad
 K_\kappa=\mathfrak d_0+
             \sum_{r\geq 1}\frac{(-\kappa)^r}{r!}h^{*r}(v)\,dv.
\end{equation}
Here $K_\kappa$ is a signed measure, and the convolution with $h_y$ is
additive.  The series is locally finite, since the support of $h^{*r}$
is contained in $[r,\infty)$.  Moreover, for every $\sigma>0$,
\[
 \int_{[0,\infty)}\e^{-\sigma v}\,d|K_\kappa|(v)
 \leq\sum_{r\geq0}\frac{\kappa^r}{r!}
       \left(\int_1^\infty\frac{\e^{-\sigma v}}{v}\,dv\right)^r
 <\infty,
\]
where $|K_\kappa|$ is the total variation measure.  Thus its Laplace
transform may be computed term by term.  If
\[
 E(s)=\int_1^\infty\frac{\e^{-sv}}{v}\,dv\qquad(\Re s>0),
\]
then
\begin{equation}\label{eq:kernel-transform}
 \widehat K_\kappa(s)=\exp\{-\kappa E(s)\}.
\end{equation}

For completeness, we justify the corresponding transform of $z_\kappa$.
Integrating \eqref{eq:z-kappa} gives
\[
 z_\kappa(v)=1-\kappa\int_1^v\frac{z_\kappa(t-1)}{t}\,dt
       \qquad(v\geq1).
\]
Consequently $M(v)=\sup_{0\leq t\leq v}|z_\kappa(t)|$ satisfies
\[
 M(v)\leq1+\kappa\int_1^v\frac{M(t)}{t}\,dt\qquad(v\geq1).
\]
The integral form of Gronwall's inequality yields $M(v)\leq v^\kappa$.
In particular, $Z(s)=\widehat z_\kappa(s)$ converges absolutely for
$\Re s>0$, as do the integrations involving $z_\kappa'$ below.
Integrating by parts and then using \eqref{eq:z-kappa} gives
\[
 -Z(s)-sZ'(s)
   =\int_0^\infty v z_\kappa'(v)\e^{-sv}\,dv
   =-\kappa\e^{-s}Z(s).
\]
For positive real $s$, the function $F(s)=sZ(s)$ therefore satisfies
\[
 F'(s)=\frac{\kappa\e^{-s}}{s}F(s),\qquad
 \lim_{s\to\infty}F(s)=1.
\]
The limit follows from $z_\kappa(v)=1$ for $0\leq v\leq1$ and the
polynomial bound just proved.  Solving this first-order equation gives
$F(s)=\exp\{-\kappa E(s)\}$.  Both sides are holomorphic on $\Re s>0$,
so the same identity holds throughout that half-plane.

Regard $dz_\kappa$ as the Stieltjes measure of the right-continuous
function extended by zero to negative arguments.  It has a unit atom
at zero and density $z_\kappa'(v)$ for $v>0$; this density also has at
most polynomial growth.  Stieltjes integration by parts gives
\[
 \widehat{dz_\kappa}(s)=sZ(s)=\widehat K_\kappa(s).
\]
Uniqueness of the Laplace transform for signed measures of finite
exponentially weighted total variation now yields $K_\kappa=dz_\kappa$.
This uniqueness statement follows, for example, by multiplying the
measures by $\e^{-\sigma v}$ and applying uniqueness of Fourier
transforms of finite measures on the line $\Re s=\sigma>0$.

Finally, $h_y$ has bounded variation on compact intervals, vanishes
for negative arguments, and is right-continuous.  Both $dh_y$ and
$dz_\kappa$ are finite signed measures on $[0,u]$.  Fubini's theorem on
the closed triangle $\{(s,t):s,t\geq0,\ s+t\leq u\}$ gives
\[
 \begin{split}
 (K_\kappa*h_y)(u)
   &=\int_{[0,u]}h_y(u-t)\,dz_\kappa(t)\\
   &=\int_{[0,u]}z_\kappa(u-v)\,dh_y(v).
 \end{split}
\]
Both atoms at zero, and any atom on the boundary of this triangle,
are retained.  Comparison with \eqref{eq:abstract-main} proves
$M_a(x,y)=\Lambda_a(x,y)$.

Finally, $y^{-1}\ll\e^{-c_3\omega(y)}$ for some $c_3>0$.  Hence
\eqref{eq:squareful-error} is absorbed in
\eqref{eq:total-replacement-error}, proving \eqref{eq:replacement}.
For $\kappa=0$, $z_0$ is the right-continuous Heaviside function, and
\eqref{eq:abstract-main} is $xh_y(u)=A(x)$.
\end{proof}

\section{Evaluation of the main term and signed moments}
\label{sec:signed}

We first supply the elementary identities for the generalized Dickman
functions that will be used both here and in Section~\ref{sec:absolute}.

\begin{lemma}[Generalized Dickman identities]\label{lem:dickman-basic}
For every $\beta>0$ and $v>0$,
\begin{equation}\label{eq:dickman-integral}
 v\rho_\beta(v)=\beta\int_{v-1}^{v}\rho_\beta(t)\,dt,
 \qquad
 0<\rho_\beta(v)\leq\frac{v^{\beta-1}}{\Gamma(\beta)}.
\end{equation}
The function $\rho_\beta$ is nonincreasing on $(0,\infty)$ when
$0<\beta\leq1$, and is strictly decreasing there when $0<\beta<1$.
Its Laplace transform is
\begin{equation}\label{eq:dickman-transform}
 \widehat\rho_\beta(s)
   =s^{-\beta}\exp\{-\beta E(s)\},\qquad
 E(s)=\int_1^\infty\frac{\e^{-sv}}{v}\,dv\quad(\Re s>0),
\end{equation}
where $s^{-\beta}$ uses the holomorphic logarithm that is real on the
positive real axis.  In particular, for $a,b>0$,
\begin{equation}\label{eq:dickman-convolution}
 (\rho_a*\rho_b)(v)=\rho_{a+b}(v)\qquad(v>0).
\end{equation}
\end{lemma}

\begin{proof}
For $0<v\leq1$, the integral identity follows directly from
\eqref{eq:rho} and the zero extension to negative arguments.
For $v>1$, the derivatives of its two sides agree almost everywhere,
since \eqref{eq:rho} gives
\[
 \frac{d}{dv}\{v\rho_\beta(v)\}
    =\beta\{\rho_\beta(v)-\rho_\beta(v-1)\}.
\]
The two sides are locally absolutely continuous on $[1,\infty)$ and
agree at $v=1$, proving the identity for every $v>0$.
If $\rho_\beta$ had a first zero $v_0>1$, continuity and its initial
positive values would make the right side of
\eqref{eq:dickman-integral} strictly positive at $v_0$, a contradiction.
Thus $\rho_\beta(v)>0$ for all $v>0$.

For $v>1$, the differential equation also yields
\[
 \frac{d}{dv}\{v^{1-\beta}\rho_\beta(v)\}
       =-\beta v^{-\beta}\rho_\beta(v-1)\leq0
       \quad\hbox{almost everywhere}.
\]
Integration from $1$ gives the upper bound in
\eqref{eq:dickman-integral}; the initial formula gives equality for
$0<v\leq1$.  If $0<\beta\leq1$, the equation for $\rho_\beta'$ and
positivity show that this derivative is nonpositive for $v>1$.
It is strictly negative there when $\beta<1$, as it is on $(0,1)$.
Local absolute continuity and continuity at $1$ give the asserted
monotonicity without requiring a derivative at that point.

The power bound proves absolute convergence of the Laplace transform
and its derivative on $\Re s>0$.  Transforming the integral identity
by Fubini's theorem gives, with $R(s)=\widehat\rho_\beta(s)$,
\[
 -R'(s)=\beta\frac{1-\e^{-s}}{s}R(s).
\]
For positive real $s\to\infty$, the initial formula on $(0,1]$ and
the power bound on $[1,\infty)$ imply $s^\beta R(s)\to1$.
The displayed first-order equation, together with this normalization
and $E'(s)=-\e^{-s}/s$, gives \eqref{eq:dickman-transform} for
positive real $s$.  Holomorphy then extends it to $\Re s>0$.
Multiplication of these transforms proves
\eqref{eq:dickman-convolution} by uniqueness of the Laplace transform.
Initially this equality holds almost everywhere.  The power bounds
at both endpoints of the convolution integral, followed by dominated
convergence after the substitution $t=vr$, show that its left side is
continuous for $v>0$; hence the equality holds for every such $v$.
\end{proof}

We use the analytic evaluation of the abstract expression
\eqref{eq:abstract-main} due to Hanrot, Tenenbaum and Wu.  The evaluation
itself does not identify a friable sum with that expression;
Lemma~\ref{lem:replacement} supplies this separate step.

\begin{lemma}[Analytic evaluation of $\Lambda_a$]
\label{lem:evaluation}
Let $\kappa,D$ be positive integers with $\kappa\leq D$, and suppose that
$|a(n)|\leq C_0d_D(n)$ for a constant $C_0>0$.  Assume that, for $\Re s>1$,
\[
 \sum_{n\geq 1}\frac{a(n)}{n^s}=\zeta(s)^\kappa G(s),
\]
where $G$ satisfies \eqref{eq:analytic-class} for some $\delta>0$.
For every fixed $\eta>0$, uniformly as $x\to\infty$ in
$(\log x)^{1+\eta}\leq y\leq x$, one has
\begin{equation}
\label{eq:evaluation}
 \Lambda_a(x,y)
   =xL^{\kappa-1}\rho_\kappa(u)
       \left\{G(1)+O\!\left(L^{-\kappa}
                        +\frac{\log(u+1)}{L}\right)\right\}.
\end{equation}
Neither multiplicativity nor nonnegativity of $a$ is required.
\end{lemma}

\begin{proof}
Use \cite[Theorem~1.1, formula~(1.24)]{HTW}, together with the
integer-parameter half-plane refinement immediately following that
theorem, in its corrected author version.  In the notation of that
paper, take the generalized zeta factor to be $\zeta^\kappa$ and the
nonnegative majorant series to be $C_0\zeta^D$.  Both zeta exponents
are positive integers; the residual factor of the majorant is
constant.  The residual factor $G$ has an absolutely convergent
Dirichlet series for $\Re s>1$: multiply the absolutely convergent
series for $\sum a(n)n^{-s}$ by that for $\zeta(s)^{-\kappa}$.
For the latter series, absolute convergence follows directly from
the convergent product $\prod_p(1+p^{-\sigma})^\kappa$ for $\sigma>1$.
Condition~\eqref{eq:analytic-class}, after decreasing
$\delta$ if necessary, verifies exactly the required holomorphy and
growth condition on the half-plane.  The leading Taylor coefficient
$a_0$ in that theorem is $G(1)$, since
\[
 \lim_{s\to 0}\frac{s^\kappa\zeta(1+s)^\kappa G(1+s)}{1+s}=G(1).
\]
The cited half-plane refinement supplies the stated range
$(\log x)^{1+\eta}\leq y\leq x$.

The class definition in \cite{HTW} assumes $G(1)\ne 0$.  If $G(1)=0$,
apply its theorem separately to $a+d_\kappa$ and to $d_\kappa$, with the
same parameter $\kappa$ in \eqref{eq:abstract-main}.  Their residual
factors are respectively $G+1$ and $1$, and both have value one at
$s=1$.  Since $d_\kappa\leq d_D$, the former coefficient sequence is
majorized by $(C_0+1)d_D$.  Subtraction and the linearity of
$\Lambda_a$ prove the additive formula \eqref{eq:evaluation} also
when $G(1)=0$.
\end{proof}

The following coarse lower bound is enough to compare the exponential
error in Lemma~\ref{lem:replacement} with the natural main-term scale.

\begin{lemma}
\label{lem:dickman-lower}
For each positive integer $\kappa$ there is $B_\kappa>0$ such that
\begin{equation}
\label{eq:rho-lower}
 \rho_\kappa(v)\geq
       \exp\{-B_\kappa v\log(2v)\}\qquad(v\geq 1).
\end{equation}
\end{lemma}

\begin{proof}
Lemma~\ref{lem:dickman-basic} supplies the integral identity,
positivity, and monotonicity of $\rho=\rho_1$ on $[0,\infty)$.
Restricting its integral identity to $[v-1,v-1/2]$ gives
\[
 \rho(v)\geq\frac{\rho(v-1/2)}{2v}\qquad(v\geq 1).
\]
Iterate until the argument lies in $(0,1]$, where $\rho=1$.  There
are at most $2v$ factors, and each denominator is at most $2v$;
therefore
\[
 \rho(v)\geq(2v)^{-2v}\qquad(v\geq 1).
\]

For integral $\kappa$, \eqref{eq:dickman-convolution} gives
$\rho_\kappa=\rho^{*\kappa}$.  For sufficiently
large $v$, restrict each of the first $\kappa-1$ convolution variables
to
\[
 \left[\frac{v}{2\kappa},\frac{v}{2\kappa}+1\right].
\]
All the arguments of $\rho$, including the last argument, then lie
in $(0,v]$.  The integration box has volume one, so monotonicity
implies $\rho_\kappa(v)\geq\rho(v)^\kappa$.  The preceding bound
proves \eqref{eq:rho-lower} for large $v$.  Positivity and continuity
of $\rho_\kappa$ on the remaining compact interval allow an increase
of $B_\kappa$ to cover every $v\geq 1$.
\end{proof}

We record explicitly the quantitative comparison underlying the
range in \eqref{eq:signed-range}.  Let $c>0$ and $B>0$ be fixed.  If
\begin{equation}
\label{eq:quantitative-signed-range}
 A_0(X\log X)^{2/3}\leq L\leq X,
\end{equation}
then, for all sufficiently large $X$,
\begin{equation}
\label{eq:quantitative-dickman-loss}
 \frac{u\log(2u)}{\sqrt L}
     =\frac{X\log(2X/L)}{L^{3/2}}
     \leq \frac{2}{A_0^{3/2}},
 \qquad
 \frac{\log X}{\sqrt L}=o(1),
\end{equation}
uniformly in this range.  Choose $A_0$ so large that
$2B/A_0^{3/2}\leq c/4$.  For every fixed $K\geq 0$, increasing the
lower threshold on $x$ then gives
\begin{equation}
\label{eq:quantitative-absorption}
 \{\log(2x)\}^K
     \exp\{-c\sqrt L+B u\log(2u)\}
       \ll \exp\{-\tfrac12c\sqrt L\}.
\end{equation}
Thus a sufficiently large fixed constant in
\eqref{eq:quantitative-signed-range} absorbs both the Dickman loss
and the logarithmic factors.

\begin{proof}[Proof of Theorem~\ref{thm:signed}]
Take $a(n)=\lambda_{\Sym^j f}(n^k)^\ell$.  Lemmas
\ref{lem:majorant} and \ref{lem:prime} permit the use of
Lemma~\ref{lem:replacement} with $\omega(t)=\sqrt{\log t}$.
Lemmas~\ref{lem:euler} and \ref{lem:unrestricted} verify the
hypotheses of Lemma~\ref{lem:evaluation}.  Its range contains
\eqref{eq:quantitative-signed-range} for all sufficiently large $x$;
for example, in that range $y\geq(\log x)^2$ eventually.

Let $c>0$ be the constant in the resulting replacement error.  Divide
that error by $xL^{\kappa-1}\rho_\kappa(u)$.  Since $\kappa\geq1$
and $L\geq1$, Lemma~\ref{lem:dickman-lower} bounds the quotient by
\[
 \ll \{\log(2x)\}^{2D}
          \exp\{-c\sqrt L+B_\kappa u\log(2u)\}.
\]
Choose $A_0$ in \eqref{eq:quantitative-signed-range} using
\eqref{eq:quantitative-absorption} with $B=B_\kappa$ and $K=2D$.
The quotient is then $O(\e^{-c\sqrt L/2})$, which is absorbed in
$O(\log(u+1)/L)$ because $u\geq1$ and $L\to\infty$.
Furthermore,
\[
 L^{-\kappa}\ll\frac{\log(u+1)}{L}.
\]
Formula~\eqref{eq:evaluation} now proves \eqref{eq:signed}, with
the constant $A$ in \eqref{eq:signed-range} chosen sufficiently
large.  The identity $C=G(1)$ and its Euler product interpretation
are supplied by Lemma~\ref{lem:euler}.

Finally, $\log(u+1)/L\to0$ uniformly in this range.  If $C\ne0$,
the additive formula therefore gives the asserted matching order
of magnitude.
\end{proof}

\begin{proof}[Proof of Theorem~\ref{thm:zero}]
For $\kappa=0$, Lemma~\ref{lem:replacement} gives
\[
 S_{j,k,\ell}(x,y;f)
    =A(x)+O\!\left(x\{\log(2x)\}^{2D}\e^{-c\sqrt L}\right).
\]
By Lemma~\ref{lem:unrestricted}, $A(x)\ll x^{1-\delta_0}$ for
some $\delta_0>0$.  Since $L\leq X$, this term is absorbed in the
displayed error, uniformly for $2\leq y\leq x$, after adjustment
of the constants.  This proves \eqref{eq:zero-global}.

The classical friable asymptotic is
\begin{equation}
\label{eq:psi}
 \Psi(x,y)=x\rho(u)
    \left\{1+O_\eta\!\left(\frac{\log(u+1)}{\log y}\right)\right\},
\end{equation}
uniformly for
$\exp\{(\log\log x)^{5/3+\eta}\}\leq y\leq x$, for every fixed
$\eta>0$; see \cite{HildebrandTenenbaum}.  In
particular, this holds throughout
\eqref{eq:quantitative-signed-range}, where its relative error
tends to zero uniformly.  Divide \eqref{eq:zero-global} by
\eqref{eq:psi} and use Lemma~\ref{lem:dickman-lower} with
$\kappa=1$.  The quotient is at most
\[
 \ll \{\log(2x)\}^{2D}
       \exp\{-c\sqrt L+B_1u\log(2u)\}.
\]
A sufficiently large choice of $A_0$, followed by
\eqref{eq:quantitative-absorption}, proves
\eqref{eq:zero-relative}, with a possibly smaller positive
constant in the exponential.
\end{proof}

\section{Nonnegative moments in a wider range}
\label{sec:nonnegative}

For nonnegative coefficients we use the mean-value theorem of
Tenenbaum and Wu.  Its real-parameter form will also be needed for
the absolute first-moment majorants in Section~\ref{sec:absolute}.

\begin{lemma}[Nonnegative friable means]
\label{lem:nonnegative}
Let $g$ be a nonnegative multiplicative function with $g(1)=1$,
$g(p)\ll1$, and $g(p^\nu)\ll(\nu+1)^B$ for some fixed $B$.
Suppose that, for fixed $\beta>0$, $c>0$, and $0<\vartheta<3/5$,
\begin{equation}
\label{eq:nonnegative-mean}
 \sum_{p\leq t}g(p)\log p
   =\beta t+O\!\left(t\e^{-c(\log t)^\vartheta}\right)
          \qquad(t\geq3).
\end{equation}
Then
\begin{equation}
\label{eq:nonnegative-constant}
 C_g=\lim_{z\to\infty}
       \prod_{p\leq z}(1-p^{-1})^\beta
          \sum_{\nu\geq0}\frac{g(p^\nu)}{p^\nu}
\end{equation}
exists and is positive.  There is a constant $A_0>0$, depending on
the fixed data in the hypotheses, such that, uniformly as
$x\to\infty$ in
\begin{equation}
\label{eq:nonnegative-range}
 \exp\{A_0(\log\log x)^{1/\vartheta}\}\leq y\leq x,
\end{equation}
one has
\begin{equation}
\label{eq:nonnegative}
 \sum_{\substack{n\leq x\\P^+(n)\leq y}}g(n)
    =C_gxL^{\beta-1}\rho_\beta(u)
       \left\{1+O\!\left(\frac{\log(u+1)}{L}+L^{-\beta}\right)\right\}.
\end{equation}
The assertion includes $0<\beta<1$.
\end{lemma}

\begin{proof}
We verify the bounded-prime clause of \cite[Theorem~2.1]{TW}.
Take
\[
 R(t)=\exp\{c_0(\log t)^\vartheta\}\qquad(t>1),
\]
with $c_0>0$ sufficiently small in terms of the constant in
\eqref{eq:nonnegative-mean}.  To check the regularity class
$\mathcal R(b;\beta)$ of that paper, put
$\phi(v)=\log R(\e^v)=c_0v^\vartheta$.  Then $\phi'$ decreases
and $v\phi'(v)$ increases.  Choose a fixed $b\in(0,1/2]$ so small
that
\[
 2^{-\vartheta}<(1-b)^2.
\]
The sufficient differential inequality in \cite[(2.6)]{TW}
follows from
\[
 \tfrac12\phi'(v/2)=2^{-\vartheta}\phi'(v)
              \leq(1-b)^2\phi'(v).
\]
The required integrability condition holds because, for $w>1$,
\[
 \int_w^\infty\frac{dt}{tR(t^b)\log t}
    =\int_{\log w}^\infty
          \frac{\exp\{-c_0b^\vartheta v^\vartheta\}}{v}\,dv<\infty.
\]
The lower growth condition in that class holds since this
exponential dominates every fixed power of $\log t$.  Its upper
growth condition holds because $\vartheta<3/5$: choose
$\delta>0$ with $\vartheta<3/5-\delta$, and then
$R(t)\ll\exp\{(\log t)^{3/5-\delta}\}$.  Thus all the stated
conditions defining $\mathcal R(b;\beta)$ are satisfied.

Fix any $\eta\in(0,1/2)$.  The prime-power hypothesis gives
\[
 \sum_p\sum_{\nu\geq2}
      \frac{g(p^\nu)\log(p^\nu)}{p^{\nu(1-\eta)}}<\infty.
\]
Indeed, the inner sum is
$O_{B,\eta}((\log p)p^{-2(1-\eta)})$, whose sum converges.
In particular the growth condition (1.8) of \cite{TW} holds.
Together with \eqref{eq:nonnegative-mean}, this places $g$ in
the class to which their theorem applies.

We spell out both its domain and its error in this specialization.
With the notation of \cite[(2.8)--(2.10)]{TW},
\[
 R_b(y)=\int_{3/2}^y\frac{dt}{tR(t^b)}\ll1,
 \qquad
 R(y^b)=\exp\{c_1L^\vartheta\},
 \qquad c_1=c_0b^\vartheta>0.
\]
Since $g(p)$ is bounded, the last term in the general error of
\cite[(2.9)]{TW} is omitted.  The remaining relative error is
\begin{equation}
\label{eq:nonnegative-intermediate-error}
 O\!\left(
       \frac{\log(u+1)}{L}+L^{-\beta}
       +u\log(u+1)\e^{-c_1L^\vartheta}
     \right).
\end{equation}
Their bounded-prime domain can be taken to be
\[
 2\leq y\leq x,\qquad
 u\leq U^*(y),\qquad
 U^*(y)=\frac{\exp\{c_1L^\vartheta\}}{1+c_1L^\vartheta},
\]
by fixing their domain parameter equal to one.

Choose $A_0$ so large that $c_1A_0^\vartheta\geq2$.
In \eqref{eq:nonnegative-range},
\[
 c_1L^\vartheta\geq c_1A_0^\vartheta\log X.
\]
The function $t\mapsto\e^t/(1+t)$ is increasing for $t\geq0$.
Consequently
\[
 U^*(y)\geq
    \frac{X^{c_1A_0^\vartheta}}
         {1+c_1A_0^\vartheta\log X}
       \geq X
\]
for all sufficiently large $X$.  Since $u=X/L\leq X$ eventually,
the required bounded-prime domain is satisfied uniformly.
Moreover, the final term in
\eqref{eq:nonnegative-intermediate-error} satisfies the exact
comparison
\[
 u\log(u+1)\e^{-c_1L^\vartheta}
    =\frac{\log(u+1)}{L}\,
          X\e^{-c_1L^\vartheta}
    \leq\frac{\log(u+1)}{L}\,
          X^{1-c_1A_0^\vartheta}.
\]
It is therefore absorbed in the first error term.  Theorem~2.1
of \cite{TW} now yields \eqref{eq:nonnegative} with its stated
constant.

For completeness, convergence and positivity of that constant can
also be seen directly.  Partial summation in
\eqref{eq:nonnegative-mean}, followed by the prime harmonic-sum
formula, shows that
\[
 \sum_p\frac{g(p)-\beta}{p}
\]
converges in increasing prime order.  The prime-power bound gives
\[
 (1-p^{-1})^\beta
       \sum_{\nu\geq0}\frac{g(p^\nu)}{p^\nu}
     =1+\frac{g(p)-\beta}{p}+O(p^{-2}).
\]
All the local factors are positive.  For sufficiently large $p$
their logarithms equal $(g(p)-\beta)/p+O(p^{-2})$; these
logarithms form a convergent series.  Keeping the finitely many
remaining positive factors outside the logarithm proves that
\eqref{eq:nonnegative-constant} exists and is strictly positive.
\end{proof}

\begin{proof}[Proof of Theorem~\ref{thm:even}]
For even $\ell$, the multiplicative function
$a(n)=\lambda_{\Sym^j f}(n^k)^\ell$ is nonnegative.
Proposition~\ref{prop:parity} gives $\kappa\geq1$.
Lemmas~\ref{lem:majorant} and \ref{lem:prime} verify the
hypotheses of Lemma~\ref{lem:nonnegative} with
$\beta=\kappa$ and $\vartheta=1/2$.  Its range is
\eqref{eq:even-range} after choosing a sufficiently large fixed
constant $A$.  Since $\kappa\geq1$ and $u\geq1$,
\[
 L^{-\kappa}\ll\frac{\log(u+1)}{L},
\]
and \eqref{eq:nonnegative} becomes \eqref{eq:even}.
The positive constant agrees with \eqref{eq:constant}, either
directly from \eqref{eq:nonnegative-constant} or from
Lemma~\ref{lem:euler}.
\end{proof}

\section{Two signed odd moments with positive constants}
\label{sec:odd}

We now establish the local positivity needed for
Corollary~\ref{cor:odd}. Lemma~\ref{lem:euler} reduces this question to
positivity of the local sums at $s=1$; in particular, no separate
nonvanishing assertion for a signed infinite product is required.

\subsection{A first moment at perfect squares}

Put $\alpha=\alpha_p$ and $b_p=\alpha^2+1+\alpha^{-2}\in[-1,3]$.
Taking the even part of \eqref{eq:local} for $j=2$ gives
\begin{equation}\label{eq:odd-square-local}
  \mathcal A_p(T):=
  \sum_{\nu\ge0}\lambda_{\Sym^2 f}(p^{2\nu})T^\nu
  =\frac{1+b_pT}{(1-T)(1-\alpha^4T)(1-\alpha^{-4}T)}.
\end{equation}
Indeed, if
\[
  B(q)=\frac{1}{(1-q)(1-\alpha^2q)(1-\alpha^{-2}q)},
\]
then $\mathcal A_p(q^2)=\{B(q)+B(-q)\}/2$. Putting the two fractions
over a common denominator gives \eqref{eq:odd-square-local}.
For real $0<T<1$, its denominator is
$(1-T)|1-\alpha^4T|^2>0$ and its numerator satisfies
$1+b_pT\ge1-T>0$. Thus $\mathcal A_p(1/p)>0$ at every prime.

The prime identity
\[
  \lambda_{\Sym^2 f}(p^2)=1+\lambda_{\Sym^4 f}(p)
\]
is also the character identity
$\Sym^2(V_2)=V_4\oplus V_0$. It gives $\kappa_{2,2,1}=1$ and
\[
  F(s)=\zeta(s)L(s,\Sym^4 f)H(s).
\]
Here, with $T=p^{-s}$,
\begin{equation}\label{eq:odd-square-correction}
  H_p(s)=(1+b_pT)(1-T)(1-\alpha^2T)(1-\alpha^{-2}T)
         =1+O(p^{-2\Re s}).
\end{equation}
The cancellation of its linear term also follows directly from
$b_p=1+\alpha^2+\alpha^{-2}$. In particular,
\[
  H_p(1)=(1+b_p/p)(1-1/p)|1-\alpha_p^2/p|^2>0.
\]
Lemma~\ref{lem:euler}, applied to the positive local sums just proved,
gives $C_{2,2,1}(f)>0$. The factorization provides the useful expression
\begin{equation}\label{eq:odd-square-constant}
  C_{2,2,1}(f)
  =L(1,\Sym^4 f)
    \prod_p(1+b_p/p)(1-1/p)|1-\alpha_p^2/p|^2>0,
\end{equation}
where the correction product converges absolutely. Applying
Theorem~\ref{thm:signed} proves \eqref{eq:odd-square}.

This is a signed example already at primes. The same local identities
give
\[
  \lambda_{\Sym^2 f}(p^2)=b_p^2-b_p<0
  \qquad\text{when }0<b_p<1.
\]
The condition $0<1+2\cos(2\theta)<1$ defines a set of positive
Sato--Tate measure. It therefore holds at a positive proportion of
primes by \cite{Thorner}. The positivity of the main constant does not
come from nonnegativity of all the coefficients.

\subsection{A cubic moment}

Set $t=\alpha_p^2$ and $b=t+1+t^{-1}\in[-1,3]$. The decomposition
\[
  V_2^{\otimes3}\simeq V_6\oplus V_4^{\oplus2}
                          \oplus V_2^{\oplus3}\oplus V_0
\]
gives $\kappa_{2,1,3}=1$. The exact local generating function is
\begin{equation}\label{eq:odd-cubic-local}
  \sum_{\nu\ge0}\lambda_{\Sym^2 f}(p^\nu)^3T^\nu
  =\frac{N_b(T)}{(1-T)\prod_{r=1}^{3}(1-t^rT)(1-t^{-r}T)},
\end{equation}
where
\begin{equation}\label{eq:odd-cubic-numerator}
  N_b(T)=1+(2b^2+b-1)(T+T^3)
                  +(b^3+3b^2-2b)T^2+T^4.
\end{equation}

Here are details of the identity. For $t\ne\pm1$, partial fractions in
the local series with parameters $t,1,t^{-1}$ give
\[
  \lambda_{\Sym^2 f}(p^\nu)
  =\frac{t^{\nu+3}+t^{-\nu}-t(t+1)}{(t-1)^2(t+1)}.
\]
Write
$(t^3Z+Z^{-1}-t(t+1))^3=\sum_{r=-3}^{3}c_r(t)Z^r$.
Direct expansion gives
\begin{align*}
  c_3(t)&=t^9,&c_2(t)&=-3t^7(t+1),\\
  c_1(t)&=3t^5(t^2+3t+1),&
  c_0(t)&=-t^3(t+1)(t^2+8t+1),\\
  c_{-1}(t)&=3t^2(t^2+3t+1),&c_{-2}(t)&=-3t(t+1),
  \qquad c_{-3}(t)=1.
\end{align*}
Cubing the preceding coefficient formula and summing the seven geometric
series yields
\[
  \sum_{\nu\ge0}\lambda_{\Sym^2 f}(p^\nu)^3T^\nu
  =\frac{1}{(t-1)^6(t+1)^3}
       \sum_{r=-3}^{3}\frac{c_r(t)}{1-t^rT}.
\]
Putting these terms over the common denominator
$\prod_{r=-3}^{3}(1-t^rT)$ and using $t+t^{-1}=b-1$ gives exactly
\eqref{eq:odd-cubic-numerator}. This proves \eqref{eq:odd-cubic-local}
for $t\ne\pm1$. Both sides depend continuously on the unit-circle
parameter $t$ when $|T|<1$: on the series side this follows from the
uniform polynomial bound on the coefficients, and on the rational side
the denominator is nonzero. Hence the identity also holds at
$t=\pm1$.

For $-1\le b\le3$,
\[
  2b^2+b-1\ge-\frac98,
  \qquad b^3+3b^2-2b\ge-\frac13.
\]
The first inequality follows by completing the square. For the second,
the cubic is nonnegative on $[-1,0]$, while on $[0,3]$ it is at least
$3b^2-2b\ge-1/3$. Consequently, for $0<T\le1/2$,
\begin{equation}\label{eq:odd-cubic-positivity}
  N_b(T)
  \ge1-\frac98\left(\frac12+\frac18\right)
         -\frac13\cdot\frac14
  =\frac{41}{192}>0.
\end{equation}
The denominator in \eqref{eq:odd-cubic-local} is
$(1-T)\prod_{r=1}^{3}|1-t^rT|^2>0$ for real $0<T<1$.
In particular, the local sum is strictly positive at $T=1/p$ for
every prime $p$. Lemma~\ref{lem:euler} proves $C_{2,1,3}(f)>0$.

For completeness, the full factorization is
\[
  F(s)=\zeta(s)L(s,\Sym^2 f)^3L(s,\Sym^4 f)^2
                  L(s,\Sym^6 f)H(s),
\]
where, with $T=p^{-s}$,
\begin{equation}\label{eq:odd-cubic-correction}
\begin{split}
  H_p(s)={}&N_b(T)(1-T)^6
       \{(1-tT)(1-t^{-1}T)\}^{5}\\
       &\hspace{18mm}\cdot
       \{(1-t^2T)(1-t^{-2}T)\}^{2}.
\end{split}
\end{equation}
Its linear term vanishes, so the correction product converges absolutely
in $\Re s>1/2$. The same positivity calculation shows $H_p(1)>0$, and
\[
  C_{2,1,3}(f)
  =L(1,\Sym^2 f)^3L(1,\Sym^4 f)^2L(1,\Sym^6 f)
      \prod_pH_p(1)>0.
\]
Theorem~\ref{thm:signed} now proves \eqref{eq:odd-cubic}, completing
the proof of Corollary~\ref{cor:odd}.

\section{Ordinary moments of CM newforms}\label{sec:cm}

In this section, $f$ is a fixed normalized holomorphic CM newform of weight
$w\ge2$, level $N$, and trivial nebentypus.  Let $K$ be its imaginary
quadratic CM field, $\mathcal O_K$ its ring of integers, and $\chi_K$ the
associated quadratic Dirichlet character; a bar denotes complex
conjugation in $K$.  For a nonzero integral ideal
$\mathfrak a\subseteq\mathcal O_K$, its norm is
$\operatorname N\mathfrak a=\#(\mathcal O_K/\mathfrak a)$.
The Hecke--Shimura construction and the characterization of CM newforms
by Hecke characters give a primitive unitary Hecke character $\psi$ of
$K$ such that
\begin{equation}\label{eq:cm-theta}
 L(s,f)=L_K(s,\psi),\qquad
 \lambda_f(n)=\sum_{\operatorname N\mathfrak a=n}\psi(\mathfrak a).
\end{equation}
Here and below, ideal sums run over nonzero integral ideals.  This is the
unitary normalization of the theta-series identity in
\cite[Section~1, Theorems~1.1 and~1.4]{Schuett}: if $\varphi$ is the
algebraic Hecke character of infinity type $w-1$ (so
$\varphi((\alpha))=\alpha^{w-1}$ for $\alpha\equiv1$ modulo its
conductor, with the fixed complex embedding of $K$), then
$\psi(\mathfrak a)=\varphi(\mathfrak a)
(\operatorname N\mathfrak a)^{-(w-1)/2}$ on ideals coprime to its conductor.
In particular, $|\psi(\mathfrak a)|=1$ on those ideals, and its infinity
type is nonzero.

For each integer $h\ge1$, let $\psi_h$ denote the primitive character
associated to the power $\psi^h$, and let $\mathfrak f_h$ be its conductor;
thus $\psi_1=\psi$.  Extend $\psi_h$ by zero on integral ideals not
coprime to $\mathfrak f_h$.  We always use its primitive Hecke $L$-function,
\begin{equation}\label{eq:hecke-definition}
 L_K(s,\psi_h)
 =\sum_{\mathfrak a\ne0}\frac{\psi_h(\mathfrak a)}
                                {(\operatorname N\mathfrak a)^s}
 =\prod_{\mathfrak p\nmid\mathfrak f_h}
       \left(1-\frac{\psi_h(\mathfrak p)}
                         {(\operatorname N\mathfrak p)^s}\right)^{-1}
 \qquad(\Re s>1),
\end{equation}
where $\mathfrak p$ ranges over nonzero prime ideals.  The conductor
$\mathfrak f_h$ can be smaller than $\mathfrak f_1$; in particular, the
zero extension of $\psi_h$ need not equal the pointwise $h$th power of
the zero extension of $\psi$.  Away from the primes above the level,
the two agree.  The infinity type of $\psi_h$ is $h$ times that of
$\psi$, so $L_K(s,\psi_h)$ is entire.

Identity \eqref{eq:cm-theta} also proves the needed local bound at every
rational prime, including those dividing the level.  There are at most
$\nu+1$ integral ideals of norm $p^\nu$: there are exactly $\nu+1$ when
$p$ splits, at most one when $p$ is inert, and one when $p$ ramifies.
Taking absolute values in \eqref{eq:cm-theta}, with the zero extension
at the conductor, therefore gives
\begin{equation}\label{eq:cm-local-bound}
 |\lambda_f(p^\nu)|\le\nu+1\qquad(\nu\ge0).
\end{equation}

All identities at rational primes in the next paragraph exclude the
finitely many primes dividing $N\operatorname{disc}(K)$.

At a split prime $p\mathcal O_K=\mathfrak p\overline{\mathfrak p}$, write
$z=\psi(\mathfrak p)$.  The two normalized Satake parameters are
$\psi(\mathfrak p)$ and $\psi(\overline{\mathfrak p})$.  Their product is
one because the nebentypus is trivial; hence
$\psi(\overline{\mathfrak p})=z^{-1}$, $|z|=1$, and
$\lambda_f(p)=z+z^{-1}$.  At an inert prime, there is no integral ideal
of norm $p$, so \eqref{eq:cm-theta} gives $\lambda_f(p)=0$.  For an
integer $h\ge1$, put
\[
 b_h(p)=
 \begin{cases}
 z^h+z^{-h},&p\text{ split in }K,\\
 0,&p\text{ inert in }K.
 \end{cases}
\]
Set $b_h(p)=0$ at the excluded primes.  At every good prime, these
are the coefficients at rational primes of the degree-two
$L$-function $L_K(s,\psi_h)$.  If $\ell$ is odd, the binomial theorem gives
\begin{equation}\label{eq:cm-odd-primes}
 \lambda_f(p)^\ell
 =\sum_{q=0}^{(\ell-1)/2}\binom{\ell}{q}b_{\ell-2q}(p).
\end{equation}
For $\ell=2r$, it gives
\begin{equation}\label{eq:cm-even-primes}
 \lambda_f(p)^{2r}
 =\kappa_{\CM,2r}\{1+\chi_K(p)\}
   +\sum_{q=0}^{r-1}\binom{2r}{q}b_{2r-2q}(p),
 \qquad \kappa_{\CM,2r}=\frac12\binom{2r}{r}.
\end{equation}
The factor $1+\chi_K(p)$ is twice the indicator of the split primes.
This accounts for the factor $1/2$ in the CM moment parameter.
For odd $\ell$, set $\kappa_{\CM,\ell}=0$.

\begin{lemma}[Effective CM prime moments]\label{lem:cm-prime}
For each fixed positive integer $\ell$, there is $c>0$ such that
\begin{equation}\label{eq:cm-prime-mean}
 \sum_{p\le t}\lambda_f(p)^\ell\log p
 =\kappa_{\CM,\ell}t
  +O_{f,\ell}\bigl(t\e^{-c\omega_{\CM}(t)}\bigr)
 \qquad(t\ge3),
\end{equation}
where
\[
 \omega_{\CM}(t)
 =\frac{(\log t)^{3/5}}{\{\log\log(3t)\}^{1/5}}.
\]
\end{lemma}

\begin{proof}
Every character $\psi_h$ occurring in
\eqref{eq:cm-odd-primes}--\eqref{eq:cm-even-primes} has nonzero infinity
type.  Coleman's zero-free region \cite[Theorem~2]{Coleman}, in the form
stated in \cite[Theorem~2.5]{Khale}, therefore gives, for each fixed such
character, a zero-free region
\begin{equation}\label{eq:cm-zero-free}
 \Re s\ge1-
 \frac{c_h}{\{\log(|\Im s|+3)\}^{2/3}
                  \{\log\log(|\Im s|+\e^{\e})\}^{1/3}}
\end{equation}
with $c_h>0$.  The conductor and the infinity type are fixed here.
The exceptional real-character case in that theorem does not occur for
$\psi_h$.  At bounded height, one obtains the same assertion after
reducing $c_h$, by nonvanishing on $\Re s=1$ and compactness.

We spell out the prime-number-theorem consequence.  Define the ideal
von Mangoldt function by
\[
 \Lambda_K(\mathfrak a)=
 \begin{cases}
 \log\operatorname N\mathfrak p,&\mathfrak a=\mathfrak p^r,
                       \ \mathfrak p\text{ prime},\ r\ge1,\\
 0,&\text{otherwise}.
 \end{cases}
\]
Write
\[
 \Psi_K(t,\psi_h)
 =\sum_{\operatorname N\mathfrak a\le t}
       \Lambda_K(\mathfrak a)\psi_h(\mathfrak a).
\]
The truncated explicit formula, at a suitable height
$T'\in[T,2T]$, gives
\begin{equation}\label{eq:cm-explicit}
 \Psi_K(t,\psi_h)
 =-\sum_{|\Im\rho|\le T'}\frac{t^\rho}{\rho}
  +O_{f,h}\left(\frac{t\log^2(tT)}{T}+\log^2t\right),
\end{equation}
where the sum runs over the nontrivial zeros of $L_K(s,\psi_h)$.
This is the usual fixed-conductor explicit formula; see, for example,
the explicit-formula and zero-counting results in \cite[Chapter~5]{IK}.
One may first use a half-integer endpoint and then move the endpoint to
$t$.  Since $K$ is quadratic, the coefficients of the ideal von Mangoldt
sum at any fixed integer norm are $O_K(\log t)$, so this change is absorbed
in the displayed error.

Let $N_h(T)$ denote the number of nontrivial zeros $\rho$ of
$L_K(s,\psi_h)$ in $0<\Re\rho<1$ with $|\Im\rho|\le T$, counted with
multiplicity.
The fixed-character bound is
\mbox{$N_h(T)\ll_{f,h}T\log(T+3)$} for $T\ge1$.  Partial summation gives
\[
 \begin{split}
 \sum_{1<|\Im\rho|\le T'}\frac1{|\rho|}
 &\le\sum_{1<|\Im\rho|\le T'}\frac1{|\Im\rho|}\\
 &\le\frac{N_h(T')}{T'}+\int_1^{T'}\frac{N_h(v)}{v^2}\,dv
 \ll_{f,h}\log^2(T+3).
 \end{split}
\]
The finitely many nontrivial zeros with $|\Im\rho|\le1$ are nonzero
and contribute only to the implied constant.  Thus the same bound
holds for $\sum_{|\Im\rho|\le T'}|\rho|^{-1}$.  Combining this with \eqref{eq:cm-zero-free} and
\eqref{eq:cm-explicit}, and decreasing the zero-free-region constant if
necessary, gives
\[
 \Psi_K(t,\psi_h)
 \ll_{f,h}
 t\exp\left\{-\frac{c_h\log t}
                   {(\log T)^{2/3}(\log\log T)^{1/3}}\right\}\log^2T
 +\frac{t\log^2(tT)}{T}+\log^2t.
\]
Choose $\log T=\omega_{\CM}(t)$ for sufficiently large $t$.
Both principal error terms are then
$O_{f,h}(t\e^{-c'_h\omega_{\CM}(t)})$, after absorbing their logarithmic
factors and choosing $c'_h>0$ sufficiently small.  Removing ideal prime
powers and prime ideals above inert rational primes costs
$O_K(\sqrt t\log^2t)$.  Hence
\begin{equation}\label{eq:cm-nontrivial-primes}
 \sum_{p\le t}b_h(p)\log p
 \ll_{f,h}t\e^{-c'_h\omega_{\CM}(t)}.
\end{equation}

For the principal Hecke character, the pole at one supplies a main term
$t$ in the same argument.  Consequently
\begin{equation}\label{eq:cm-split-primes}
 2\sum_{\substack{p\le t\\p\text{ split in }K}}\log p
 =t+O_K\bigl(t\e^{-c_K\omega_{\CM}(t)}\bigr).
\end{equation}
A possible exceptional real zero for the fixed field contributes
$O_K(t^{\beta_0})$, with fixed $\beta_0<1$, and is absorbed in this error.
Substituting \eqref{eq:cm-nontrivial-primes} and
\eqref{eq:cm-split-primes} into the binomial identities proves
\eqref{eq:cm-prime-mean}.  The omitted bad primes affect only the implied
constant.
\end{proof}

The same prime identities give, initially for $\Re s>1$,
\begin{equation}\label{eq:cm-factorization}
 \sum_{n\ge1}\frac{\lambda_f(n)^\ell}{n^s}
 =\zeta(s)^{\kappa_{\CM,\ell}}
   L(s,\chi_K)^{\kappa_{\CM,\ell}}
   \prod_{q=0}^{\lfloor(\ell-1)/2\rfloor}
      L_K(s,\psi_{\ell-2q})^{\binom{\ell}{q}}
   H_\ell(s).
\end{equation}
Indeed, after dividing the arithmetic Euler factor by the indicated
$L$-factors, its linear term in $p^{-s}$ vanishes at every good prime.
The bound \eqref{eq:cm-local-bound} then shows that the correction
factor is $1+O_{\ell,\eta}(p^{-2\Re s})$ on
$\Re s\ge1/2+\eta$.  Thus the Euler product for $H_\ell$ converges
absolutely and locally uniformly on $\Re s>1/2$, and it is bounded on
every closed half-plane $\Re s\ge1/2+\eta$.
The finitely many bad-prime correction factors in this product use the
actual local factors of $f$ and the primitive local factors of
$L_K(s,\psi_h)$.  They are holomorphic for $\Re s>0$: by
\eqref{eq:cm-local-bound}, the arithmetic local series has coefficients
of polynomial growth, and the inverses of the indicated local
$L$-factors are finite products of $1-\xi p^{-rs}$ with $|\xi|\le1$
and $r\in\{1,2\}$.

The elementary dimension argument of Lemma~\ref{lem:majorant}, applied
to the ordinary coefficient sequence, gives
\begin{equation}\label{eq:cm-majorant}
 |\lambda_f(n)|^\ell\le d_{2^\ell}(n).
\end{equation}
Equivalently, the local inequality follows by mapping multisets of
$\nu$ binary vectors in $\{0,1\}^\ell$ onto their coordinate sums:
every element of $\{0,\ldots,\nu\}^\ell$ occurs.  Thus
$(\nu+1)^\ell\le\binom{2^\ell+\nu-1}{\nu}$, and multiplicativity proves
\eqref{eq:cm-majorant}.

The total degree of the displayed $L$-factors in
\eqref{eq:cm-factorization} is $2^\ell$.  For odd $\ell$, this follows
from $2\sum_{q=0}^{(\ell-1)/2}\binom{\ell}{q}=2^\ell$.
For $\ell=2r$, it follows from
\[
 2\kappa_{\CM,2r}
 +2\sum_{q=0}^{r-1}\binom{2r}{q}=2^{2r}.
\]
When $\ell$ is odd, the zeta factor is absent and every displayed Hecke
$L$-function is entire.  Convexity for these fixed degree-two
$L$-functions and \eqref{eq:cm-majorant} therefore give, by the same
truncated Perron argument as in Lemma~\ref{lem:unrestricted},
\begin{equation}\label{eq:cm-unrestricted}
 \sum_{n\le t}\lambda_f(n)^\ell
 \ll_{f,\ell}t^{1-\delta_\ell}
 \qquad(\ell\text{ odd})
\end{equation}
for some $\delta_\ell>0$.  More explicitly, put $d=2^\ell$, truncate at
$T=t^{2/(d+4)}$, and shift to $\Re s=1/2+\eta$ with a sufficiently
small fixed $\eta>0$.  Convexity bounds the vertical integral by
$t^{1/2+\eta}T^{d(1/2-\eta)/2+\varepsilon}$, while the horizontal and
truncation terms are $O_\varepsilon(t^{1+O(\varepsilon)}/T)$.
Both have a fixed power saving when $\eta$ and $\varepsilon$ are
sufficiently small.  There is no residue.

\begin{proof}[Proof of Theorem~\ref{thm:cm}]
For an even moment, Lemma~\ref{lem:cm-prime} supplies the hypothesis
\eqref{eq:nonnegative-mean} for every fixed $\vartheta<3/5$.
Given $\varepsilon>0$, choose $\vartheta$ sufficiently close to $3/5$
that $1/\vartheta<5/3+\varepsilon$.
The strict inequality absorbs the fixed range constant in
Lemma~\ref{lem:nonnegative} for sufficiently large $x$.
Since $\kappa_{\CM,2r}\ge1$, the additional error
$L^{-\kappa_{\CM,2r}}$ is absorbed by $\log(u+1)/L$.
This proves \eqref{eq:cm-even}, including positivity and the stated
Euler-product expression for its constant.  The actual bad-prime local
sums are retained in that expression.

For an odd moment, apply Lemma~\ref{lem:replacement} with
$D=2^\ell$, $\kappa=0$, and $\omega=\omega_{\CM}$.
Its abstract main term is the unrestricted sum in
\eqref{eq:cm-unrestricted}.  Since
$\omega_{\CM}(y)\ll\omega_{\CM}(x)=o(\log x)$ for $2\le y\le x$,
the power-saving unrestricted error is absorbed, after a possible
decrease of the exponential constant.  Thus, for some $c_0>0$,
\begin{equation}\label{eq:cm-global-proof}
 |T_\ell(x,y;f)|
 \ll_{f,\ell}
 x\{\log(2x)\}^{2^{\ell+1}}
 \e^{-c_0\omega_{\CM}(y)}.
\end{equation}
This proves \eqref{eq:cm-odd-global}.

To obtain the relative estimate in the stated stronger range, set
$X=\log x$ and $W=\log X$, and suppose
\begin{equation}\label{eq:cm-balanced-range}
 A X^{5/8}W^{3/4}\le L\le X,
\end{equation}
where the fixed constant $A\ge1$ will be chosen sufficiently large.
For all sufficiently large $x$, uniformly in this range,
\begin{align}
 \omega_{\CM}(y)
 &\ge c_2 A^{3/5}X^{3/8}W^{1/4},\label{eq:cm-omega-comparison}\\
 u\log(2u)
 &\le 2A^{-1}X^{3/8}W^{1/4}.\label{eq:cm-u-comparison}
\end{align}
Here $c_2>0$ is absolute: the first estimate uses
$\log\log(3y)\ll W$, and the second uses $\log(2u)\le2W$.
In particular,
\[
 u\log(2u)\ll A^{-8/5}\omega_{\CM}(y),
 \qquad
 W=o\bigl(\omega_{\CM}(y)\bigr).
\]
The range \eqref{eq:cm-balanced-range} lies inside the classical friable
range of \eqref{eq:psi}.  By that estimate and
Lemma~\ref{lem:dickman-lower}, there is a fixed $B_1>0$ such that
\[
 \frac{|T_\ell(x,y;f)|}{\Psi(x,y)}
 \ll_{f,\ell}
 \exp\left\{
 2^{\ell+1}\log\log(2x)+B_1u\log(2u)
       -c_0\omega_{\CM}(y)\right\}.
\]
Choose $A$ sufficiently large that the middle positive term is at most
$c_0\omega_{\CM}(y)/4$.  For sufficiently large $x$, the first positive
term is also at most $c_0\omega_{\CM}(y)/4$.
The last display is consequently
$O_{f,\ell}(\e^{-c_0\omega_{\CM}(y)/2})$, proving
\eqref{eq:cm-odd-relative} in \eqref{eq:cm-odd-range}.
\end{proof}

For ordinary even moments, the first few prime means are
\[
\begin{array}{c|rrrr}
 \ell&2&4&6&8\\ \hline
 \text{non-CM }\kappa_{1,1,\ell}&1&2&5&14\\
 \text{CM }\kappa_{\CM,\ell}&1&3&10&35
\end{array}
\]
The difference between these parameters is why the two hypotheses on
$f$ must be distinguished.

\section{Absolute majorants for first moments}\label{sec:absolute}

We prove Corollary~\ref{cor:first}.  For a level-one form, the prime
coefficients in its two families agree:
\[
 \lambda_{\Sym^m f}(p)=\lambda_f(p^m)=U_m(\cos\theta_p).
\]
Character orthogonality gives
\[
 \int_0^\pi U_m(\cos\theta)^2\,d\mu_{\ST}(\theta)=1,
 \qquad d\mu_{\ST}(\theta)=\frac2\pi\sin^2\theta\,d\theta.
\]
Strict Cauchy--Schwarz therefore yields $\mu_m<1$, since
$|U_m(\cos\theta)|$ is continuous and nonconstant on the support of
this probability measure.

The mean also has the closed form
\begin{equation}\label{eq:mu-explicit}
 \mu_m=\frac{4(m+1)}{\pi m(m+2)}
            \cot\frac{\pi}{2(m+1)}.
\end{equation}
To verify it, put $n=m+1\ge2$.  The defining integral becomes
\[
 \mu_m=\frac2\pi\int_0^\pi|\sin(n\theta)|\sin\theta\,d\theta.
\]
On the interval $[r\pi/n,(r+1)\pi/n]$, replace
$|\sin(n\theta)|$ by $(-1)^r\sin(n\theta)$.
Integrating the product of the two sine functions and summing the
endpoint values gives
\[
 \int_0^\pi|\sin(n\theta)|\sin\theta\,d\theta
 =\frac{2n}{n^2-1}\sum_{r=1}^{n-1}\sin\frac{r\pi}{n}
 =\frac{2n}{n^2-1}\cot\frac{\pi}{2n},
\]
which proves \eqref{eq:mu-explicit}.  In particular,
$\mu_1=8/(3\pi)$.

For the CM case, we derive the limiting prime distribution directly
from Lemma~\ref{lem:cm-prime}.  We use the Vinogradov--Korobov form of
the ordinary prime number theorem:
\begin{equation}\label{eq:vk-pnt}
 \sum_{p\le t}\log p
 =t+O\bigl(t\e^{-c\omega_{\CM}(t)}\bigr)\qquad(t\ge3),
\end{equation}
for an absolute $c>0$.  For example, this follows from
\cite[Theorem~1.4]{JohnstonYang}, after absorbing its power of $\log t$
by decreasing $c$.
For a continuous function $F$ on $[-2,2]$, set
\[
 \mathcal M_{\CM}(F)=\frac12F(0)
       +\frac1{2\pi}\int_0^\pi F(2\cos\theta)\,d\theta.
\]
This is a probability mean.  Its odd moments vanish, and its $2r$th
moment, for $r\ge1$, is $\tfrac12\binom{2r}{r}$.  Indeed, symmetry
about $\pi/2$ gives the odd moments, and expansion of
$(\e^{i\theta}+\e^{-i\theta})^{2r}$ shows that only the constant
Fourier term survives integration over $[0,\pi]$.
Consequently, Lemma~\ref{lem:cm-prime} and \eqref{eq:vk-pnt} imply
\begin{equation}\label{eq:cm-continuous-mean}
 \frac1t\sum_{p\le t}F(\lambda_f(p))\log p
       \longrightarrow\mathcal M_{\CM}(F).
\end{equation}
To justify the assertion first for a polynomial, expand it in monomials
and use those moment identities; its constant term uses \eqref{eq:vk-pnt}.
For general $F$, choose a polynomial $Q$ with
$\sup_{[-2,2]}|F-Q|<\delta$.  Since
$|\lambda_f(p)|\le2$ by \eqref{eq:cm-local-bound}, replacing $F$ by $Q$
in the left side of \eqref{eq:cm-continuous-mean} changes it by at most
$\delta t^{-1}\sum_{p\le t}\log p=\delta+o(1)$, and changes its
proposed limit by at most $\delta$.  Letting $\delta\to0$ proves the
claim.  Taking $F(t)=|t|$ gives
\begin{equation}\label{eq:mu-cm}
 \mu_{\CM}=\mathcal M_{\CM}(|\cdot|)
 =\frac1{2\pi}\int_0^\pi|2\cos\theta|\,d\theta=\frac2\pi.
\end{equation}
In particular, the mass $1/2$ at zero in this mean has been accounted
for explicitly.

Fix $m\ge1$ and $0<\eta<1-\mu_m$.  Polynomial approximation on
$[-2,2]$ gives a fixed even polynomial $P$ such that
\begin{equation}\label{eq:absolute-polynomial}
 |U_m(t/2)|\le P(t)\le |U_m(t/2)|+\eta
 \qquad(-2\le t\le2).
\end{equation}
For example, first approximate the continuous even function
$|U_m(t/2)|$ uniformly within $\eta/3$, take the even part of the
approximating polynomial, and add $2\eta/3$.
Let
\[
 \beta=\int_0^\pi P(2\cos\theta)\,d\mu_{\ST}(\theta).
\]
Then $0<\mu_m\le\beta\le\mu_m+\eta<1$.
Expanding $P(t)$ in the basis $U_h(t/2)$ and using
Lemma~\ref{lem:prime} for the finitely many nonconstant basis terms,
together with \eqref{eq:vk-pnt} for the constant term, gives
\begin{equation}\label{eq:absolute-prime-mean}
 \sum_{p\le t}P(\lambda_f(p))\log p
 =\beta t+O_{f,m,\eta}(t\e^{-c\sqrt{\log t}})
\end{equation}
for some $c>0$.

Define a nonnegative multiplicative function $g$ by
\[
 g(1)=1,\qquad
 g(p)=P(\lambda_f(p)),\qquad
 g(p^\nu)=(\nu+1)^{m+1}\quad(\nu\ge2).
\]
At primes, \eqref{eq:absolute-polynomial} majorizes both coefficient
families.  At higher prime powers, the local Euler factors give
\[
 |\lambda_{\Sym^m f}(p^\nu)|
 \le\binom{\nu+m}{m}\le(\nu+1)^m,
 \qquad
 |\lambda_f(p^{m\nu})|\le m\nu+1\le(\nu+1)^m.
\]
The binomial bound counts vectors of $m$ nonnegative exponents with
sum at most $\nu$ inside $\{0,\ldots,\nu\}^m$; the second bound follows
from the binomial expansion of $(\nu+1)^m$.
Multiplicativity consequently gives
\begin{equation}\label{eq:absolute-domination}
 |\lambda_{\Sym^m f}(n)|\le g(n),
 \qquad |\lambda_f(n^m)|\le g(n).
\end{equation}

Lemma~\ref{lem:nonnegative}, with prime mean $\beta$ and
$\vartheta=1/2$, applies to this function.  Hence, for a sufficiently
large fixed $A=A(f,m,\eta)>0$, uniformly for
\begin{equation}\label{eq:absolute-noncm-range}
 \exp\{A(\log\log x)^2\}\le y\le x,
\end{equation}
we have
\begin{equation}\label{eq:absolute-g-bound}
 \sum_{\substack{n\le x\\P^+(n)\le y}}g(n)
 \ll_{f,m,\eta}xL^{\beta-1}\rho_\beta(u).
\end{equation}
This invocation uses the real-parameter version of the lemma; there is
no requirement that $\beta\ge1$.

The comparison with the ordinary Dickman function follows from
Lemma~\ref{lem:dickman-basic}: for $0<\beta<1$, $\rho_\beta$ is positive
and decreasing on $(0,\infty)$, and
$\rho_\beta*\rho_{1-\beta}=\rho$.
In particular, for $u\ge1$,
\begin{align}
 \rho(u)
 &=\int_0^u\rho_\beta(u-v)\rho_{1-\beta}(v)\,dv\notag\\
 &\ge\rho_\beta(u)
        \int_0^1\frac{v^{-\beta}}{\Gamma(1-\beta)}\,dv
 =\frac{\rho_\beta(u)}{\Gamma(2-\beta)}.
 \label{eq:rho-absolute-comparison}
\end{align}
The singularity at $v=0$ is integrable because $\beta<1$.
Combining \eqref{eq:absolute-g-bound},
\eqref{eq:rho-absolute-comparison}, and \eqref{eq:psi}, we obtain
\[
 \sum_{\substack{n\le x\\P^+(n)\le y}}g(n)
 \ll_{f,m,\eta}\Psi(x,y)L^{\beta-1}
 \ll_{f,m,\eta}\Psi(x,y)L^{\mu_m-1+\eta}.
\]
Together with \eqref{eq:absolute-domination}, this proves
\eqref{eq:first-noncm}, with the sum of the two absolute values bounded
by twice the displayed majorant.

For the CM assertion, choose a fixed even polynomial $P$ majorizing
$|t|$ within $\eta$ on $[-2,2]$, where
$0<\eta<1-2/\pi$.  Its prime mean is
\[
 \beta=\frac12P(0)
       +\frac1{2\pi}\int_0^\pi P(2\cos\theta)\,d\theta,
 \qquad \frac2\pi\le\beta\le\frac2\pi+\eta<1.
\]
The identities \eqref{eq:cm-even-primes} and
Lemma~\ref{lem:cm-prime}, applied to the finitely many even monomials
of $P$, give an effective prime mean with error
$O_{f,\eta}(t\e^{-c\omega_{\CM}(t)})$; the constant monomial uses the
quantitative prime number theorem \eqref{eq:vk-pnt}.
Set $g(p)=P(\lambda_f(p))$ outside the finite set of primes omitted in
Section~\ref{sec:cm}, use $g(p)=|\lambda_f(p)|$ at those omitted primes,
and set $g(p^\nu)=\nu+1$ for all $\nu\ge2$.
Then $|\lambda_f(n)|\le g(n)$, and the finite changes at primes do not
alter the prime mean.  Given $\varepsilon>0$, choose
$\vartheta<3/5$ so close to $3/5$ that
$1/\vartheta<5/3+\varepsilon$.
Lemma~\ref{lem:nonnegative} applies throughout
\eqref{eq:cm-even-range} for sufficiently large $x$.
The comparison \eqref{eq:rho-absolute-comparison} now proves
\eqref{eq:first-cm}.

There is also an explicit polynomial choice for $m=1$:
\begin{equation}\label{eq:explicit-absolute-polynomial}
 P_0(t)=\frac49+\frac{11}{18}t^2-\frac1{18}t^4.
\end{equation}
For $|t|\le2$,
\[
 P_0(t)-|t|
 =(|t|-1)^2\left\{\frac12-\frac{(|t|+1)^2}{18}\right\}\ge0,
 \qquad P_0(t)\le2.
\]
The second inequality follows, for instance, by differentiating
$P_0$ on $[0,2]$ and using its evenness.
The second and fourth prime moments give means $17/18$ in the
non-CM case and $8/9$ in the CM case.  Thus the same argument, without
polynomial approximation, gives
\[
 \sum_{\substack{n\le x\\P^+(n)\le y}}|\lambda_f(n)|
 \ll_f\Psi(x,y)L^{-1/18}
\]
for level-one forms in \eqref{eq:absolute-noncm-range} with a sufficiently
large $A=A(f)$, and
\[
 \sum_{\substack{n\le x\\P^+(n)\le y}}|\lambda_f(n)|
 \ll_{f,\varepsilon}\Psi(x,y)L^{-1/9}
\]
for CM forms in \eqref{eq:cm-even-range}.
These are estimates for absolute moments and do not assert optimality
for the corresponding signed sums.

\section{Prime sums and the limits of the cancellation estimates}
\label{sec:obstruction}

When the prime mean is positive, Theorem~\ref{thm:signed} determines
the upper order whenever its leading constant is nonzero.
Corollary~\ref{cor:odd} gives two signed odd moments with positive
constants.  For zero prime mean, obtaining an optimal signed bound is
a different problem.  The following exact identity explains one
obstruction to strong uniform power savings.

\begin{proposition}[A prime-sum obstruction]\label{prop:obstruction}
For every arithmetic function $a$, put $A(x)=\sum_{n\le x}a(n)$.
If $x\ge2$ and $x/2\le y\le x$, then
\begin{equation}\label{eq:prime-obstruction}
 \sum_{\substack{n\le x\\P^+(n)\le y}}a(n)
 =A(x)-\sum_{y<p\le x}a(p).
\end{equation}
In particular, suppose $a(n)=\lambda_f(n)$ for a level-one form as
above, $A(x)\ll x^\alpha$, and the left side of
\eqref{eq:prime-obstruction} is $O(x^\theta)$ uniformly for
$x/2\le y\le x$.  If $\gamma=\max\{\alpha,\theta\}>0$, then
\begin{equation}\label{eq:obstruction-prime-bounds}
 \sum_{p\le x}\lambda_f(p)\ll x^\gamma,
 \qquad
 \sum_{p\le x}\lambda_f(p)\log p\ll x^\gamma\log(2x).
\end{equation}
Consequently $L(s,f)$ has no zeros in
$\Re s>\max\{\gamma,1/2\}$.
\end{proposition}

\begin{proof}
If $n\le x$ has a prime divisor $p>y\ge x/2$, then $n=p$.
This proves \eqref{eq:prime-obstruction}.
Taking $y=x/2$ shows that
$\sum_{x/2<p\le x}\lambda_f(p)\ll x^\gamma$.
Summing over the dyadic intervals
$(x/2^{j+1},x/2^j]$ gives the first estimate in
\eqref{eq:obstruction-prime-bounds}, because $\gamma>0$.
Partial summation gives the second estimate.

That weighted prime bound implies, by partial summation and local
uniform convergence, that
\[
 P_f(s)=\sum_p\frac{\lambda_f(p)\log p}{p^s}
\]
is holomorphic for $\Re s>\gamma$.
For $\Re s>1$, the local Euler factors give
\[
 -\frac{L'}{L}(s,f)
 =P_f(s)+
   \sum_p\sum_{\nu\ge2}
      \frac{(\alpha_p^\nu+\alpha_p^{-\nu})\log p}{p^{\nu s}}.
\]
Since $|\alpha_p|=1$, the double series converges absolutely and
locally uniformly for $\Re s>1/2$.
The right side therefore supplies a holomorphic continuation of
$-L'/L$ to $\Re s>\max\{\gamma,1/2\}$.
The logarithmic derivative of an entire $L$-function has a pole at
every zero, so no zero can lie in this half-plane.
\end{proof}

\end{document}